\documentclass[12pt]{amsart}

\usepackage[utf8]{inputenc} % Allow non-ASCII input
\usepackage[T1]{fontenc}    % Allow Unicode output
\usepackage{lmodern}        % Fix [ and ] in the T1 encoding
\usepackage{amssymb}        % Blackboard bold and other 
\usepackage{amsthm}         % Custom Theorems
\usepackage{mathtools}      % Extend AMS packages
\usepackage{mathrsfs}
\usepackage{physics}        % Use shorthand for derivatives
\usepackage{asymptote}      % In-Code drawings
\usepackage{subcaption}     % Side-by-side figures
\usepackage{xcolor}         % Color Links
\usepackage[
    hyperfootnotes=false      %% Disable footnote patch warning
]{hyperref}                 % Hyperlinks in PDF
\usepackage{cleveref}       % Named references
\usepackage{microtype}      % Small typographic adjustments
\usepackage{booktabs}       % Nicely-styled tables
\usepackage{tabularx}       % Full-width tables
\usepackage[
    maxbibnames=99,           %% Print all names in bibliography 
    maxcitenames=2,           %% ... but not in the text
    backend=biber             %% use Biber instead of bibtex
]{biblatex}                 % Better bibliography
\usepackage{filecontents}   % Include the .bib file inside the .tex file
\usepackage{mathptmx}
\usepackage[margin=1.5in]{geometry}
\usepackage{float}
\usepackage{layouts}

\DeclareFieldFormat[article]{pages}{#1}
\DeclareFieldFormat[article]{volume}{\mkbibbold{#1}}
\makeatletter
\def\l@section{\@tocline{1}{0pt}{0pc}{5pc}{}}
\def\l@subsection{\@tocline{2}{0pt}{2.5pc}{5pc}{}}
\makeatother
\definecolor{linkblue}{HTML}{00356B}
\definecolor{linkgold}{HTML}{DA9100}
\definecolor{linkred}{RGB}{159, 29, 53}
\hypersetup{
        breaklinks = true, 
	colorlinks = true,
	linkcolor = linkred,
	citecolor = linkgold,
	urlcolor = linkblue    
}

\theoremstyle{plain}
\newtheorem{theorem}{Theorem}[section]
\crefname{theorem}{Theorem}{Theorems}

\crefname{conjecture}{Conjecture}{Conjectures}
\newtheorem{proposition}[theorem]{Proposition} 
\crefname{proposition}{Proposition}{Propositions}
 
\crefname{corollary}{Corollary}{Corollaries}
 
\crefname{lemma}{Lemma}{Lemmas}
\crefname{ineq}{inequality}{inequalities}
\newtheorem{remark}[theorem]{Remark}
\crefname{remark}{Remark}{Remarks}

\theoremstyle{definition}

\crefname{example}{Example}{Examples}
\newtheorem{definition}[theorem]{Definition}

\crefname{appendix}{Appendix}{Appendices}
\crefname{section}{Section}{Sections}
\crefname{figure}{Figure}{Figures}
\crefname{table}{Table}{Tables}

\newcommand{\field}[1]{\mathbf{#1}}
\newcommand{\R}{\field{R}}

\title[]{Numerically Computed Double, Triple, and Quadruple Planar Bubbles for Density $r^p$}

\author[]{Marcus Collins}

\address{Harvard-Westlake School, 700 N Faring Rd, Los Angeles, CA 90077}
\email{marcuscollins2026@gmail.com}

\date{\today}

\begin{document}

% Compiling bibliography for arXiv

% ASYMPTOTE Global Definitions
\begin{asydef}
import geometry;
import patterns;

// Drawing Colors
pen backgroundcolor = white;

picture dot(real Hx=1mm, real Hy=1mm, pen p=currentpen)
{
  picture tiling;
  path g=(0,0)--(Hx,Hy);
  draw(tiling,g,invisible);
  dot(tiling, (0,0), p+linewidth(1));
  return tiling;
}

// Add fill patterns
add("hatch", hatch(1mm, currentpen));
add("dot",dot(1mm, currentpen));

pair bubble(real radius1, real radius2, pen drawingpen = currentpen, int lab = 0) {

	// Draws the appropriate double bubble for bubbles of volume
	// pi * radius1^2 and pi * radius2^2, respectively
	// Calculate the length of the segment connecting the centers of the
	// circles, from the law of cosines. Note that the angle opposite this
	// length formed by the radii is 60 degrees.
	real length = sqrt(radius1^2 + radius2^2 - radius1 * radius2);

	pair origin = (0,0);
	path circle1 = circle(origin, radius1);
	path circle2 = circle((length,0), radius2);

	// Check if the circles have equal curvature. If so, create a straight-line
	// path between their intersection points. Otherwise, draw the bulge as
	// an arc between the intersection points with curvature equal to the 
	// difference of the two circles.
	pair[] i_points = intersectionpoints(circle1, circle2);
	path bulge;
	
	if (radius1 == radius2) {
		// If the radii are equal, the "bulge" will actually be a circle 
		// through infinity, which mathematically is just a line but
		// computationally is an error. In this case, just draw the segment
		// connecting the two points.
		bulge = i_points[0] -- i_points[1];
	} else {
		
		// Calculate the radius of the middle circular arc
		real bulge_radius = 1 / abs(1/radius1 - 1/radius2);

		// Calculate the center of the middle circular arc
		// 	(1) Draw two circles of bulge_radius centered
		//	    at the two singularities. They will intersect
		//	    one another at the two possible centers.
		//	(2) Check which center is correct by comparing
		//	    the original radii.
		path big_circ_1 = circle(i_points[0], bulge_radius);
		path big_circ_2 = circle(i_points[1], bulge_radius);
		pair[] big_i_points = intersectionpoints(big_circ_1, big_circ_2);
		pair big_center;
		if (radius1 > radius2) {
			big_center = big_i_points[1];
		} else {
			big_center = big_i_points[0];
		}
		
		// Calculate the angles from the horizontal between the center of
		// the bulge and the two singularities. Draw an arc at the center
		// using those angles. Since arcsin(...) has a limited range,
		// we must add 180 degrees to theta if the center is on the
		// right, along with reversing the path direction.
		real height = arclength(i_points[0] -- i_points[1]) / 2.0;
		real theta = aSin(height/bulge_radius);
		
		if (radius1 > radius2) {
			bulge = arc(big_center, bulge_radius, -theta + 180, theta + 180);
		} else {
			bulge = arc(big_center, bulge_radius, theta, -theta);
		}
	}

	draw(bulge, drawingpen);

	// Calculate the parameterized time when the og circles intersect
	real[][] i_times = intersections(circle1, circle2);
	real time_1_a = i_times[0][0]; // Time on circle1 of first intersection 
	real time_1_b = i_times[1][0]; // Time on circle1 of second intersection 
	real time_2_a = i_times[0][1]; // Time on circle2 of first intersection 
	real time_2_b = i_times[1][1]; // Time on circle2 of second intersection 

	// The second path is weird to account for the fact that paths always 
	// start at time 0, and we can't run the path in reverse. This then 
	// calculates how long it takes to walk the whole path, and then breaks 
	// the arc up into subpaths.
	// path bubble1 = subpath(circle1, time_1_a, time_1_b) -- bulge -- cycle;
	path bubble1 = subpath(circle1, time_1_a, time_1_b) -- reverse(bulge) -- cycle;
	path bubble2 = subpath(circle2, 0, time_2_a) -- bulge -- subpath(circle2, time_2_b, arctime(circle2, arclength(circle2))) -- cycle;

	filldraw(bubble1, pattern("hatch"), drawingpen);
	filldraw(bubble2, pattern("dot"), drawingpen);

	if (lab == 1) {
		label(Label("$\bm{\Omega_1}$",Fill(backgroundcolor)),origin);
		label(Label("$\bm{\Omega_2}$",Fill(backgroundcolor)),(length, 0));
	}
	
	return i_points[0];
}
\end{asydef}
% END ASYMPTOTE

\begin{abstract}
Using Brakke's Evolver, we numerically verify conjectured optimal planar double bubbles for density $r^p$ and provide conjectures for triple and quadruple bubbles. 
\end{abstract}

\maketitle
% \tableofcontents

% \begin{figure}
%     \centering
%     \includegraphics[width=0.7\textwidth]{Diagrams/conjecture.pdf}
%     \caption{A standard double bubble with vertex at the origin is the 
%     conjectured double bubble in $\R^n$ with density $r^p$ \cite[Fig. 1]{Morgan2021}.}
%     \label{fig:Conjectured-Double-Bubble}
% \end{figure}

\section{Introduction}

\begin{figure}[h!]
    \begin{subfigure}[t]{0.49\linewidth}
        \centering
        \includegraphics[width=0.8\linewidth]{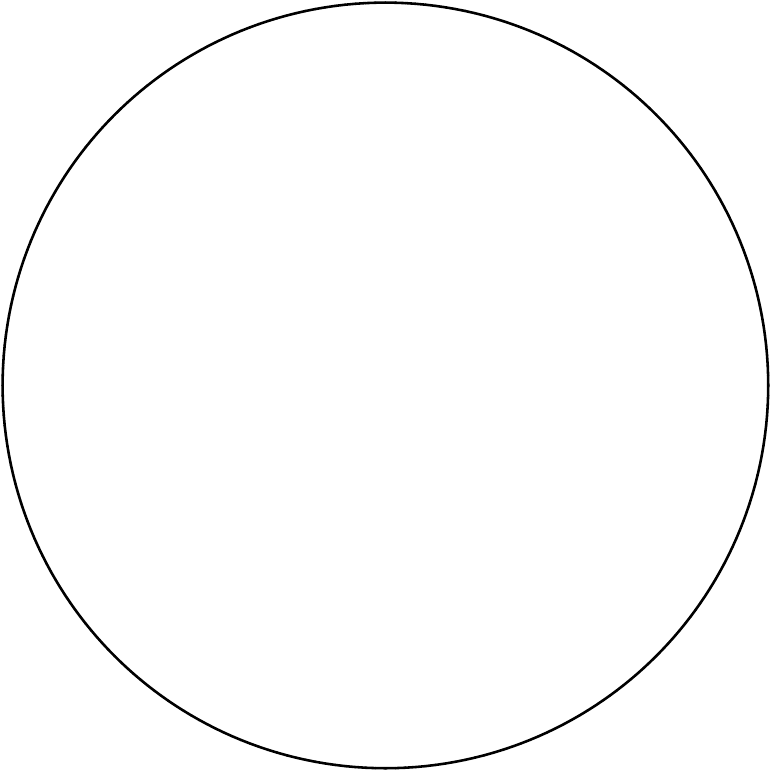}
    \end{subfigure}
    \hfill 
    \begin{subfigure}[t]{0.49\linewidth}
        \centering
        \includegraphics[width=\linewidth]{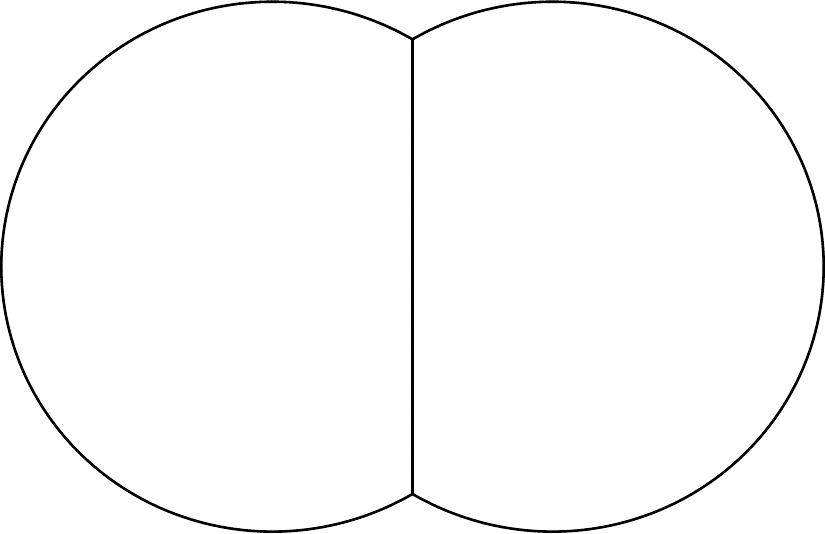}
        \begin{minipage}{.1cm}
        \vfill 
            \end{minipage} 
    \end{subfigure}
    \newline
    \begin{subfigure}[t]{0.49\linewidth}
        \centering
        \includegraphics[width=0.8\linewidth]{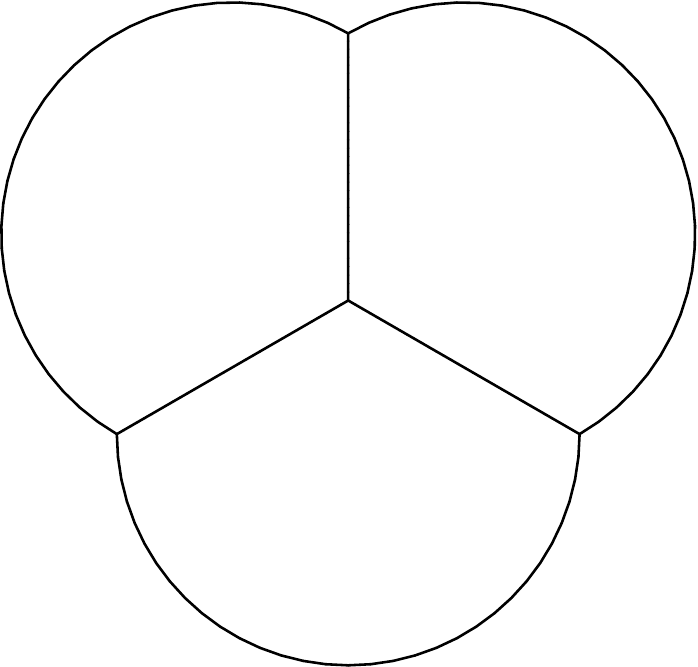}
    \end{subfigure}
    \hfill
    \begin{subfigure}[t]{0.49\textwidth}
        \centering
        \includegraphics[width=0.85\linewidth]{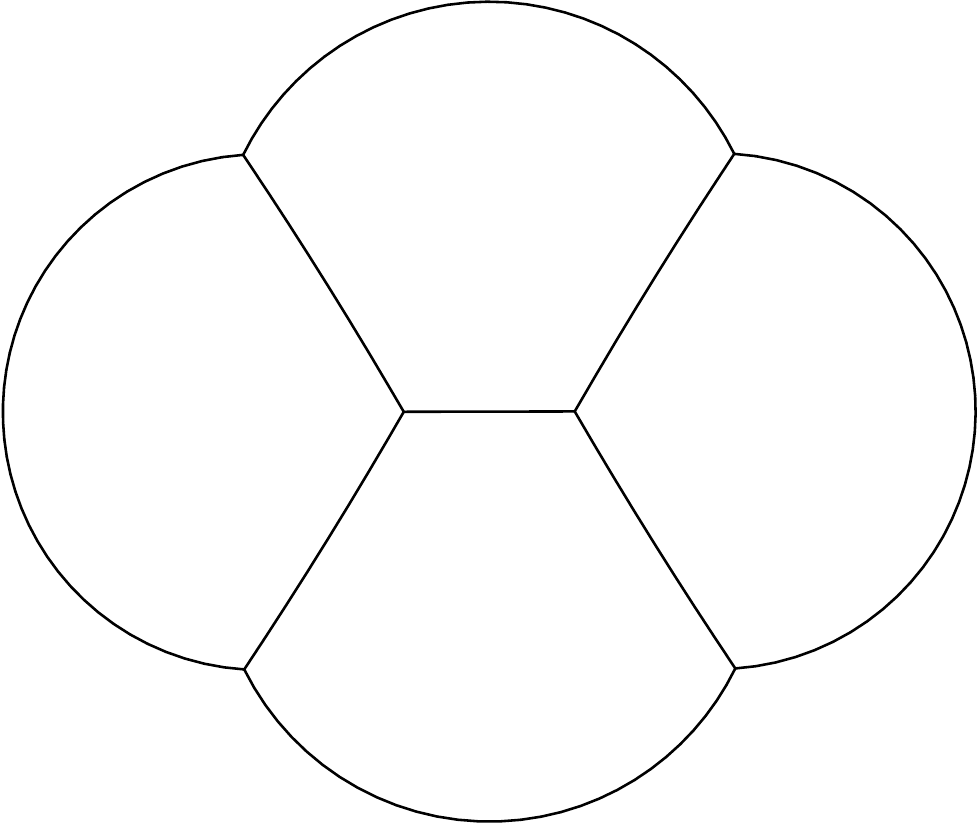}
        \begin{minipage}{.1cm}
        \vfill
        \end{minipage}
    \end{subfigure}
    \caption{The optimal Euclidean single, double \cite{foisy1993}, triple \cite{wichiramala}, and quadruple \cite{Paolini2020} bubbles with equal areas.}
    \label{fig:euclidean-bubbles}
\end{figure}

The isoperimetric problem is one of the oldest in mathematics. 
It asks for the least-perimeter way to enclose given volume. For a single volume in Euclidean space (with uniform density) of any dimension, the well-known solution is any sphere. With 
density $r^p$, \textcite{G14} found that the solution for a single volume is a sphere \emph{through} the origin. 
For \emph{two} volumes in
Euclidean space, \textcite{Reichardt2007}
showed that the standard double bubble, consisting of three spherical caps meeting along a sphere in threes at $120^\circ$ angles as in \cref{fig:euclidean-bubbles}, provides an 
isoperimetric cluster. \textcite{Morgan2021} conjecture that the isoperimetric cluster for two volumes in $\R^n$ with density $r^p$ for $p > 0$ is the same Euclidean standard double bubble with a vertex at the origin, as in \cref{fig:double-bubbles}, and show that it is better for example than putting the center at the origin. But it is not even known whether each region and the whole cluster are connected. As for the triple bubble, the minimizer in the plane with density $r^p$ cannot 
just be the Euclidean minimizer \cite{wichiramala} with central vertex at the 
origin, because the outer arcs do not have constant generalized curvature.  

\textcite{Morgan2021} prove existence, boundedness, and regularity: a planar isoperimetric cluster consists of constant generalized-curvature curves meeting in threes at $120^\circ$ (see our \cref{existence-boundedness-regularity}).

In this paper we numerically compute double, triple, and quadruple bubbles in the plane with density $r^p$ for various $p > 0$, using Brakke's Evolver \cite{Brakke2013}. Some videos are available on \href{https://drive.google.com/drive/folders/1QXaZbVwWbwMMKjbyg5eRBk6jf_D4Vi1P?usp=sharing}{Google Drive}.  \cref{doubleprop} supports the conjecture of Hirsch et al. \cite{Morgan2021} that the optimal double bubble is the Euclidean one with one vertex at the origin (\cref{fig:double-bubbles}). \cref{tripleprop} indicates that the optimal triple bubble resembles the Euclidean one (\cref{fig:euclidean-bubbles}) with one vertex at the origin, but as $p$ increases one circular arc from the origin shrinks so that all arcs pass near the origin and remain approximately circular, as in \cref{fig:triple-bubbles}. (A constant-generalized-curvature curve is a circle if and only if it passes through the origin (\cref{genconstantcurvatureremark})). \cref{quadprop} indicates that the optimal quadruple bubble resembles the Euclidean one (\cref{fig:euclidean-bubbles}) with a vertex at the origin for small $p$ (\cref{fig:quadruple-bubbles-p<1}), but as $p$ increases to $1$, one circular arc from the origin shrinks to a point and thereafter four arcs meet at the origin (\cref{fig:quad-bubbles}). This not does violate regularity because the density vanishes at the origin.
\subsection*{History} Examination of isoperimetric regions in the plane with 
density $r^p$ began in \citeyear{G06} when \textcite{G06} showed that
the isoperimetric solution for a single area in the plane with density $r^p$ is 
a convex set containing the origin. It was something of a surprise when 
\textcite{dahlberg2010} proved that the solution is a circle
through the origin.
In \citeyear{G14} \textcite{G14} extended this result to higher dimensions. In \citeyear{china19} \textcite{china19} studied the $1$-dimensional case,
showing that the best single bubble is an interval with one endpoint at the 
origin and that the best double bubble is a pair of adjacent intervals which meet at the origin. Ross \cite{Ross2022} showed that in $\R^1$ multiple bubbles start with the two smallest meeting at the origin and the rest in increasing order alternating side to side.

\subsection*{Acknowledgements}

This paper was written under the guidance of Frank Morgan remotely, following separate work done at the MathPath summer camp. The author would like to thank Professor Morgan for his advice and support. Additionally, the author would like to thank Professor Ken Brakke for generous help with his Evolver.

\section{Definitions}

\begin{definition}[Density Function]
Given a smooth Riemannian manifold $M$, a \emph{density} on $M$ is a positive 
continuous function (perhaps vanishing at isolated points) that weights each point $p$ in $M$ with a certain mass 
$f(p)$. Given a region $\Omega \subset M$ with piecewise smooth boundary, the weighted volume (or area) and boundary measure or perimeter of $\Omega$ are given by
\[
V(\Omega) = \int_{\Omega} f \dd{V_0} \quad\text{and}\quad P(\Omega) = 
\int_{\partial \Omega} f \dd{P_0},
 \]
where $\dd{V_0}$ and $\dd{P_0}$ denote Euclidean volume and perimeter. We may also refer to the perimeter of $\Omega$ as the perimeter of its boundary.
\end{definition}

\begin{definition}[Isoperimetric Region]
A region $\Omega \subset M$ is 
\emph{isoperimetric} if it has the smallest weighted perimeter of all regions 
with the same weighted volume. The boundary of an isoperimetric region is also 
called isoperimetric.
\end{definition}

We can generalize the idea of an isoperimetric region by considering two or more
volumes.
\begin{definition}[Isoperimetric Cluster]
An isoperimetric cluster is a set of disjoint open regions~$\Omega_i$ of 
volume~$V_i$ such that 
the perimeter of the union of the boundaries is minimized.
\end{definition}

For example, in the plane with density 1, optimal clusters are known for one area, two areas (Foisy et al. \cite{foisy1993}), three areas (Wichiramala \cite{wichiramala}), and four equal areas (Paolini et al. \cite{Paolini2020}), as in \cref{fig:euclidean-bubbles}. Note that for density $r^p$, scalings of minimizers are minimizers, because scaling up by a factor of $\lambda$ scales perimeter by $\lambda^p$ and area by $\lambda^{2p}$. \\

The following proposition summarizes existence, boundedness, and regularity of isoperimetric clusters in $\R^n$ with density, proved by Hirsch et al. \cite{Morgan2021} following such results for single bubbles (Rosales et al. \cite[Thm.~2.5]{rosales}) and Morgan and Pratelli \cite[Thm.~5.9]{Morgan2013}).

\begin{proposition}[Existence, Boundedness, and Regularity]\label{existence-boundedness-regularity}
    Consider $\R^n$ with radial non-decreasing $C^1$ density f such that $f(r) \to \infty$ as $r \to \infty.$ An isoperimetric cluster that encloses and separates given volumes exists (Hirsch et al.  \cite[Thm.~2.5]{Morgan2021} ) and is bounded (Hirsch et al. \cite[Prop.~2.6]{Morgan2021}). In $\R^2$, the cluster consists of smooth curves with constant generalized curvature meeting in threes at $120^\circ$ except possibly at points where the density vanishes (Hirsch et al. \cite[Thm.~2.8]{Morgan2021}). The $C^1$ hypothesis may be allowed to fail e.g. at isolated points. 
\end{proposition}
(Unlike Hirsch et al. \cite{Morgan2021}, by definition we allow a density to vanish at isolated points.)

\begin{definition}[Generalized Curvature]\label{def:gen_curve}
In $\R^2$ with density $f$, the 
generalized curvature $\kappa_f$ of a curve with inward-pointing unit normal $N$ is given by the formula
\[
\kappa_f = \kappa_0 - \pdv{\log f}{N},
\]
where $\kappa_0$ is the (unweighted) geodesic curvature. This comes from the first variation formula, so that generalized curvature has the interpretation as minus the perimeter cost $\dd{P}\!/\!\dd{A}$ of moving area across the curve, and constant generalized mean curvature is the equilibrium condition $\dd{P} = 0$ if $\dd{A} = 0$ (see \textcite[Sect.\ 3]{rosales}).
\end{definition}

\begin{remark}
    In $\R^2$ with density $r^p$, a circular arc has constant generalized curvature if and only if the circle passes through the origin  \cite[Rem.~2.9]{Morgan2021}.
    \label{genconstantcurvatureremark}
\end{remark} 
\section{Multiple Bubbles in \texorpdfstring{$\R^2$}{r^2}  with density \texorpdfstring{$r^p$}{rp}}

%Building on work by \textcite{foisy1993} and \textcite{dahlberg2010}, 

With Brakke's Evolver \cite{Brakke2013}, \cref{doubleprop} supports the conjecture of Hirsch et al. \cite{Morgan2021} that the optimal double bubble in the plane with density $r^p$ is the standard double bubble. \cref{tripleprop} provides a conjecture on the form of the triple bubble, and \cref{quadprop} provides a conjecture on the form of the quadruple bubble. 

\begin{figure}[H]
    \centering
    \begin{subfigure}[c]{0.3\textwidth}
        \centering
        \includegraphics[width=\linewidth]{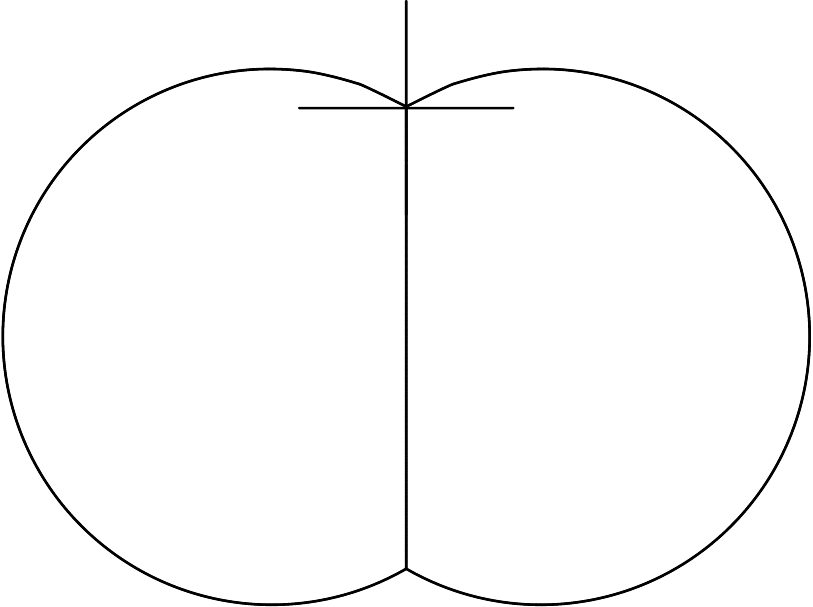}
        \caption{Equal areas}
        \label{fig:equal-double-bubble}
    \end{subfigure}
    \hfill 
    \begin{subfigure}[c]{0.3\textwidth}
        \centering
        \includegraphics[width=\linewidth]{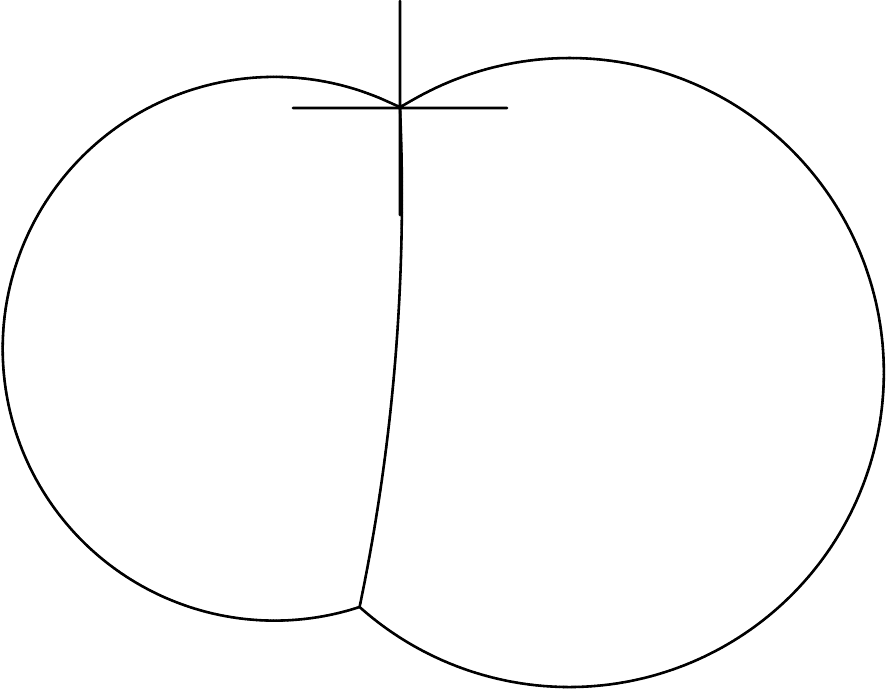}
        \caption{Areas $10$ and $20$}
        \label{fig:unequal-double-bubble}
    \end{subfigure}
    \hfill 
    \begin{subfigure}[c]{0.3\textwidth}
        \centering
        \includegraphics[width=1.2\linewidth]{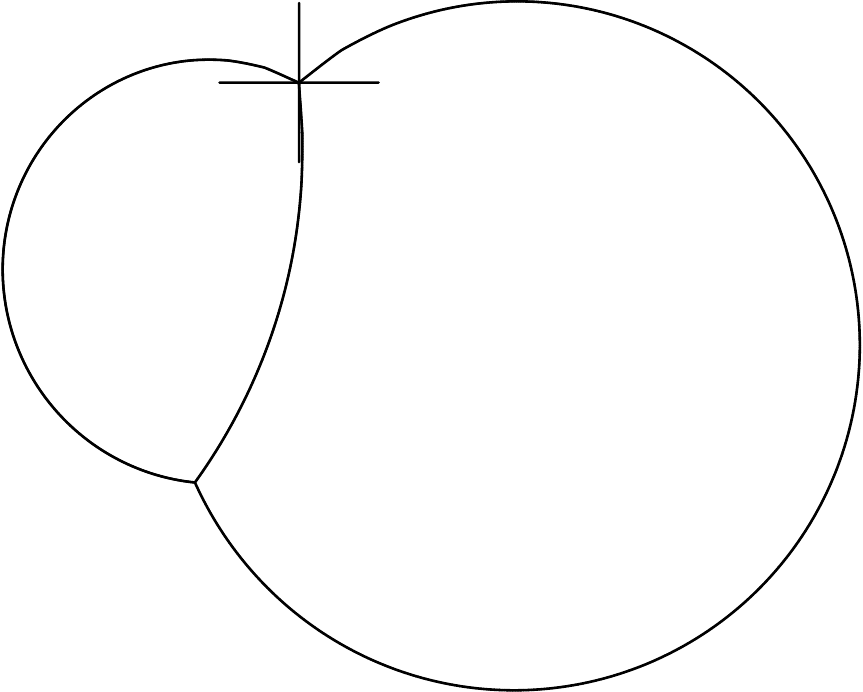}
        \caption{Areas $10$ and $100$}
        \label{fig:double-bubble-3}
        \begin{minipage}{.1cm}
        \vfill
        \end{minipage}
    \end{subfigure}
    \newline
    \begin{subfigure}[c]{0.3\textwidth}
        \centering
        \includegraphics[width=\linewidth]{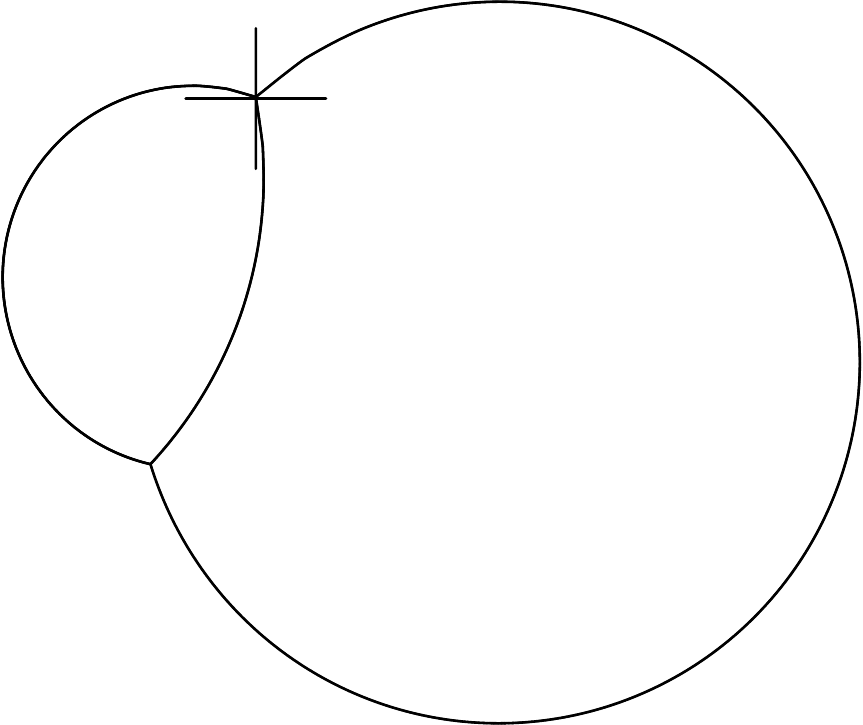}
        \caption{Areas $10$ and $200$}
        \label{fig:double-bubble-4}
        \begin{minipage}{.1cm}
        \vfill
        \end{minipage}
    \end{subfigure}
    \hfill
    \begin{subfigure}[c]{0.3\textwidth}
        \centering
        \includegraphics[width=\linewidth]{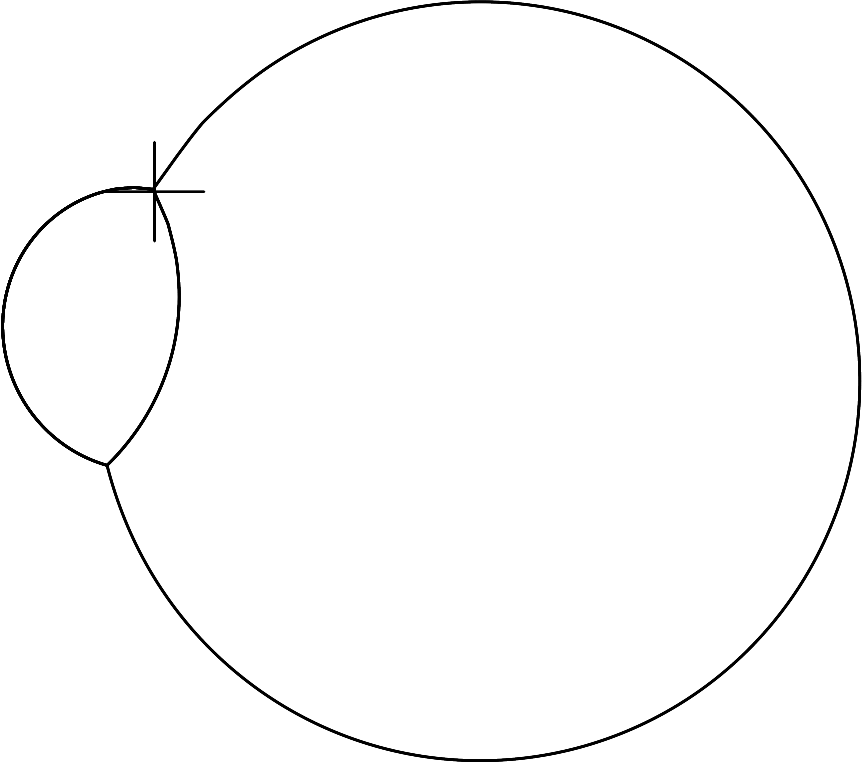}
        \caption{Areas $10$ and $1000$}
        \label{fig:double-bubble-5}
    \end{subfigure}
    \hfill
    \begin{subfigure}[c]{0.3\textwidth}
        \centering
        \includegraphics[width=\linewidth]{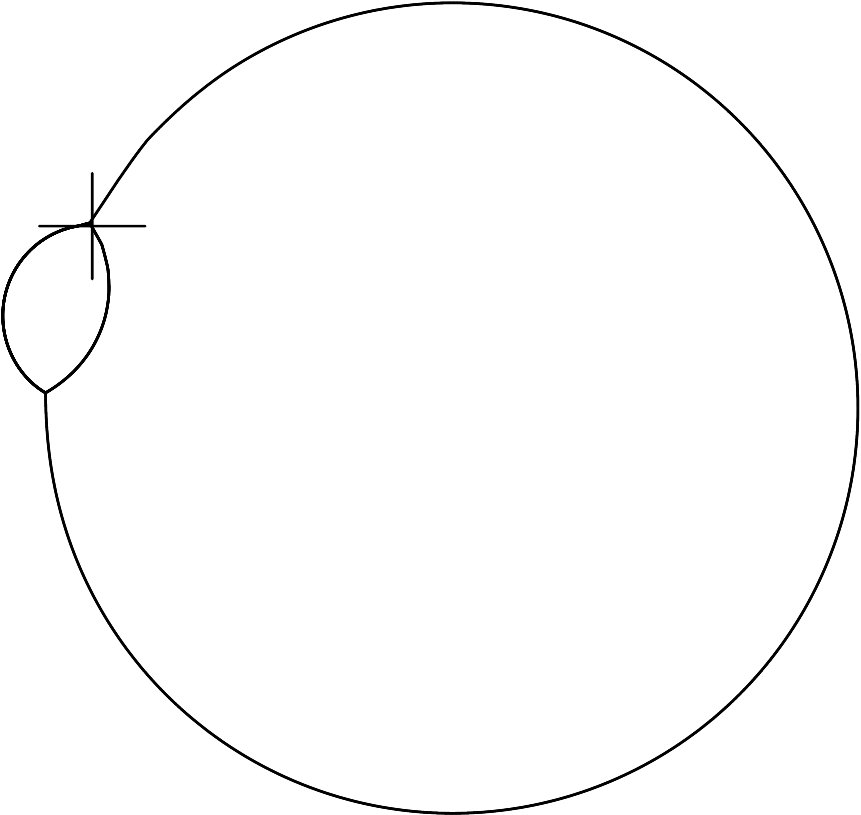}
        \caption{Areas $1$ and $1000$}
        \label{fig:double-bubble-6}
        \begin{minipage}{.1cm}
        \vfill
        \end{minipage}
    \end{subfigure}
    \caption{Computations in Brakke's Evolver \cite{Brakke2013} in $\R^2$ with density $r^2$ support the conjecture that the optimal double bubble is the standard double bubble with one vertex at the origin (marked here by a plus). Densities $r^5, r^3$, and $r^{0.5}$ are apparently identical.}
    \label{fig:double-bubbles}
\end{figure}

\begin{proposition}[Double Bubble] \label{doubleprop} 
Computations on Brakke's Evolver \cite{Brakke2013} support the conjecture \cite{Morgan2021} that the optimal planar double bubble for density $r^p$ consists of a standard double bubble with one vertex at the origin, as shown in \cref{fig:double-bubbles}.
\end{proposition}

\begin{figure}[H]
    \centering
    \begin{subfigure}[c]{0.3\textwidth}
        \centering
        \includegraphics[width=\linewidth]{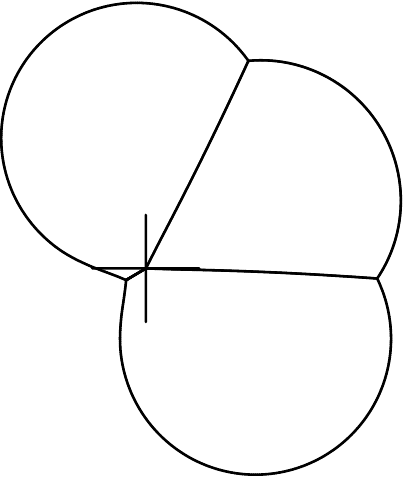}
        \caption{Equal areas of $10$}
        \label{fig:triple-bubble-1}
    \end{subfigure}
    \hfill
    \begin{subfigure}[c]{0.3\textwidth}
        \centering
        \includegraphics[width=\linewidth]{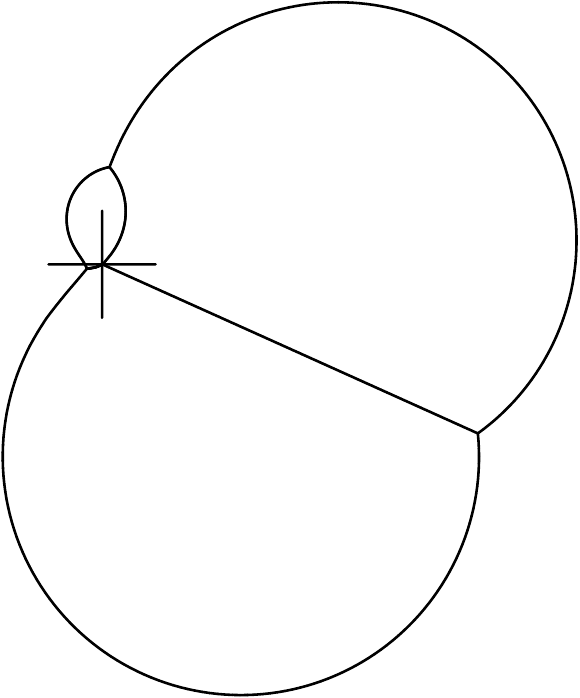}
        \caption{Areas $0.1$, $100$, and $100$}
        \label{fig:triple-bubble-5}
    \end{subfigure}
    \hfill 
    \begin{subfigure}[c]{0.3\textwidth}
        \centering
        \includegraphics[width=\linewidth]{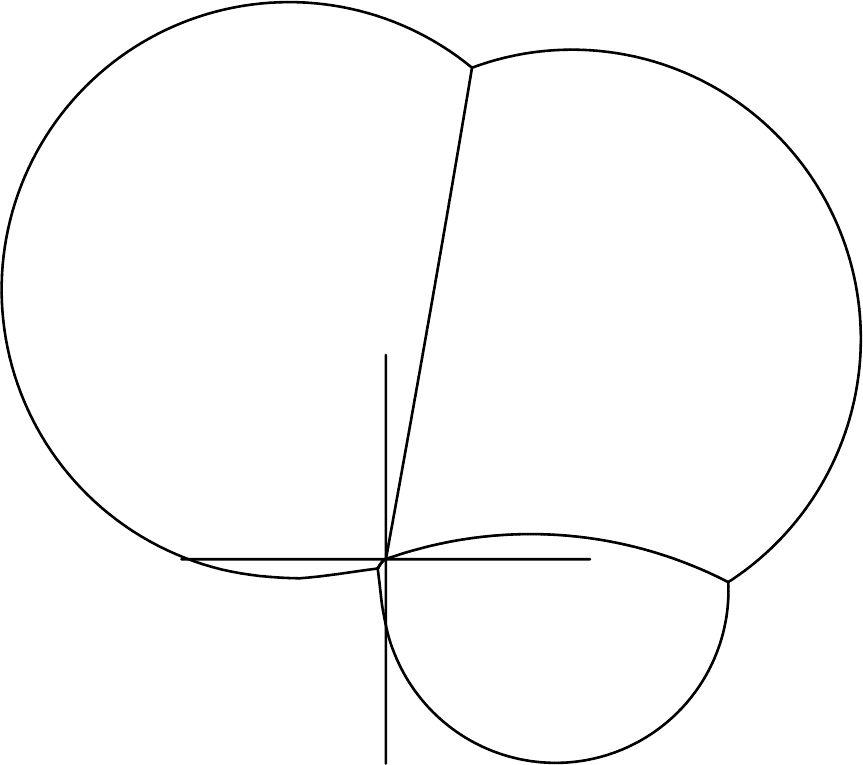}
        \caption{Areas $1$, $1$, and $0.1$}
        \label{fig:triple-bubble-3}
        \begin{minipage}{1.cm}
        \vfill
        \end{minipage}
    \end{subfigure}
    \newline 
    % \begin{subfigure}[c]{0.3\textwidth}
    %     \centering
    %     \includegraphics[height=1.257in]{Diagrams/3-2-1.png}
    %     \caption{Areas $3$, $2$, and $1$}
    %     \label{fig:triple-bubble-4}
    % \end{subfigure}
    \begin{subfigure}[c]{0.45\textwidth}
        \centering
        \includegraphics[width=\linewidth]{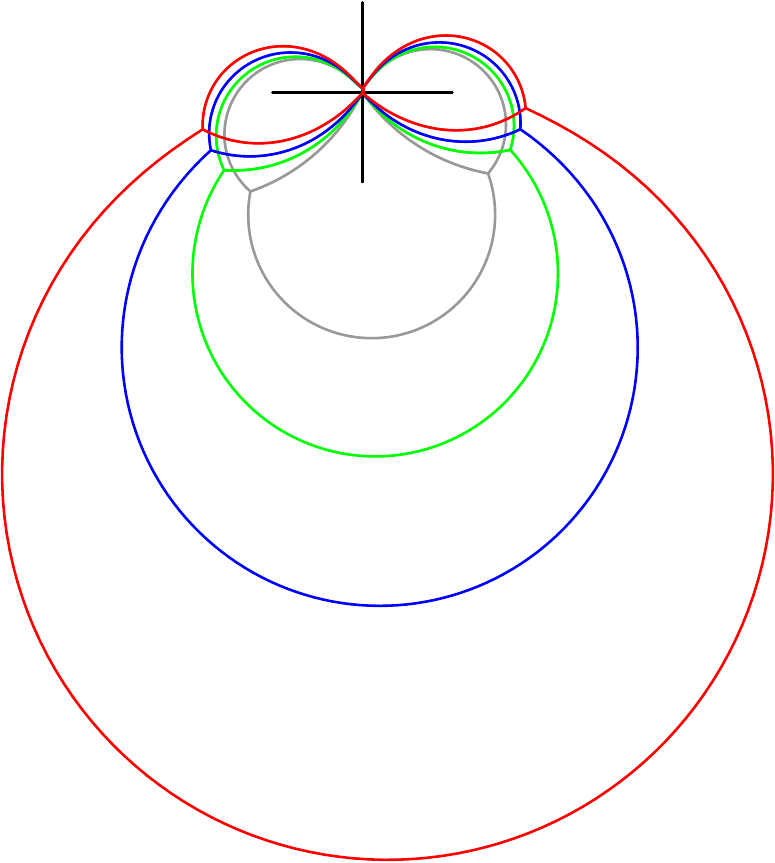}
        \caption{Areas $0.1$, $0.1$, and $1$ (grey), $0.1$, $0.1$, and $5$ (green), $0.1$, $0.1$ and $20$ (blue), and $0.1$, $0.1$ and $100$ (red) overlayed.}
        \label{fig:triple-bubble-overlayed}
    \end{subfigure}
    \hfill 
    \begin{subfigure}[c]{0.45\textwidth}
        \centering
        \includegraphics[width=\linewidth]{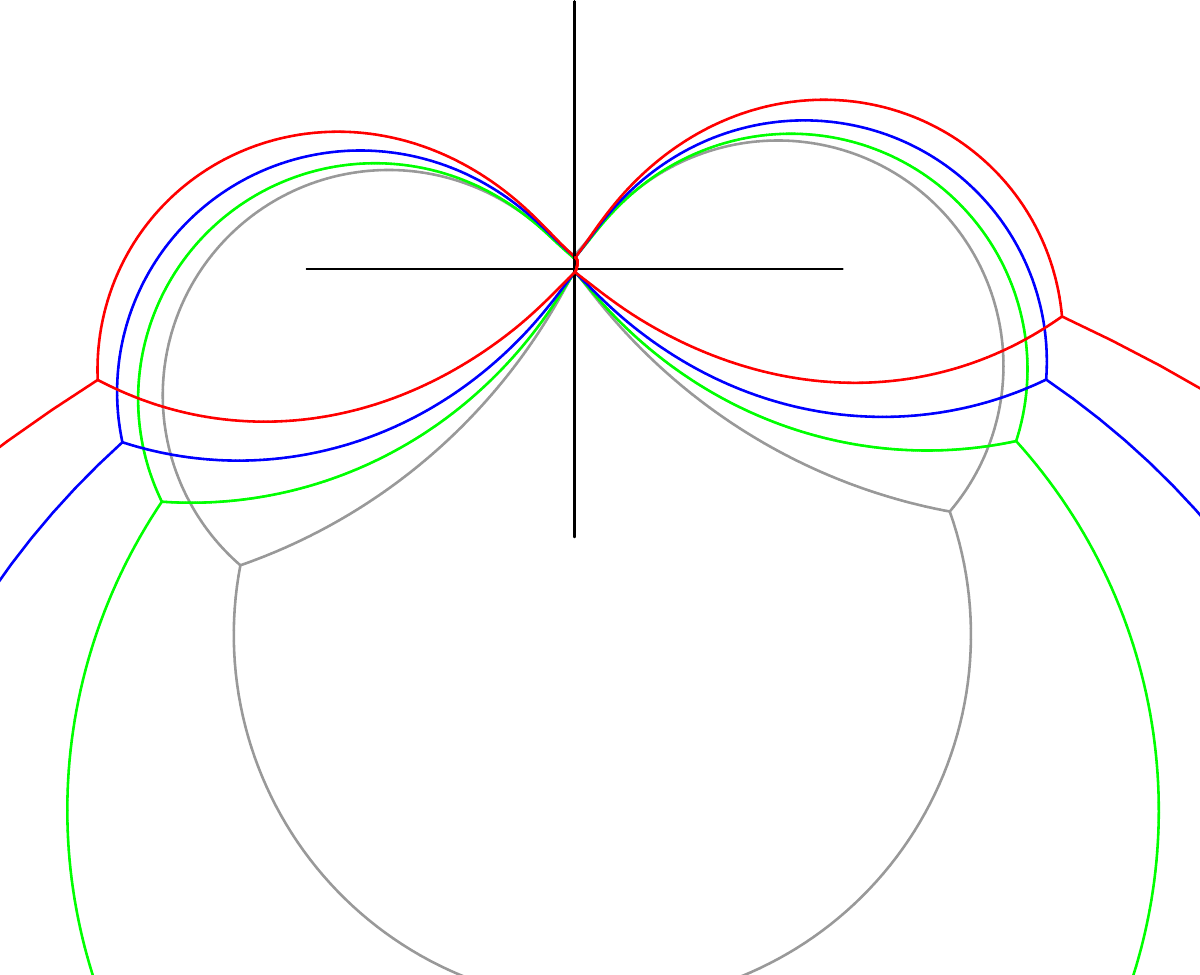}
        \caption{\cref{fig:triple-bubble-overlayed} zoomed in on the origin}
        \label{origin-overlay}
        \begin{minipage}{.1cm}
        \vfill
        \end{minipage}
    \end{subfigure} 
    \newline
    \begin{subfigure}[c]{0.2\textwidth}
        \centering
        \includegraphics[width=\linewidth]{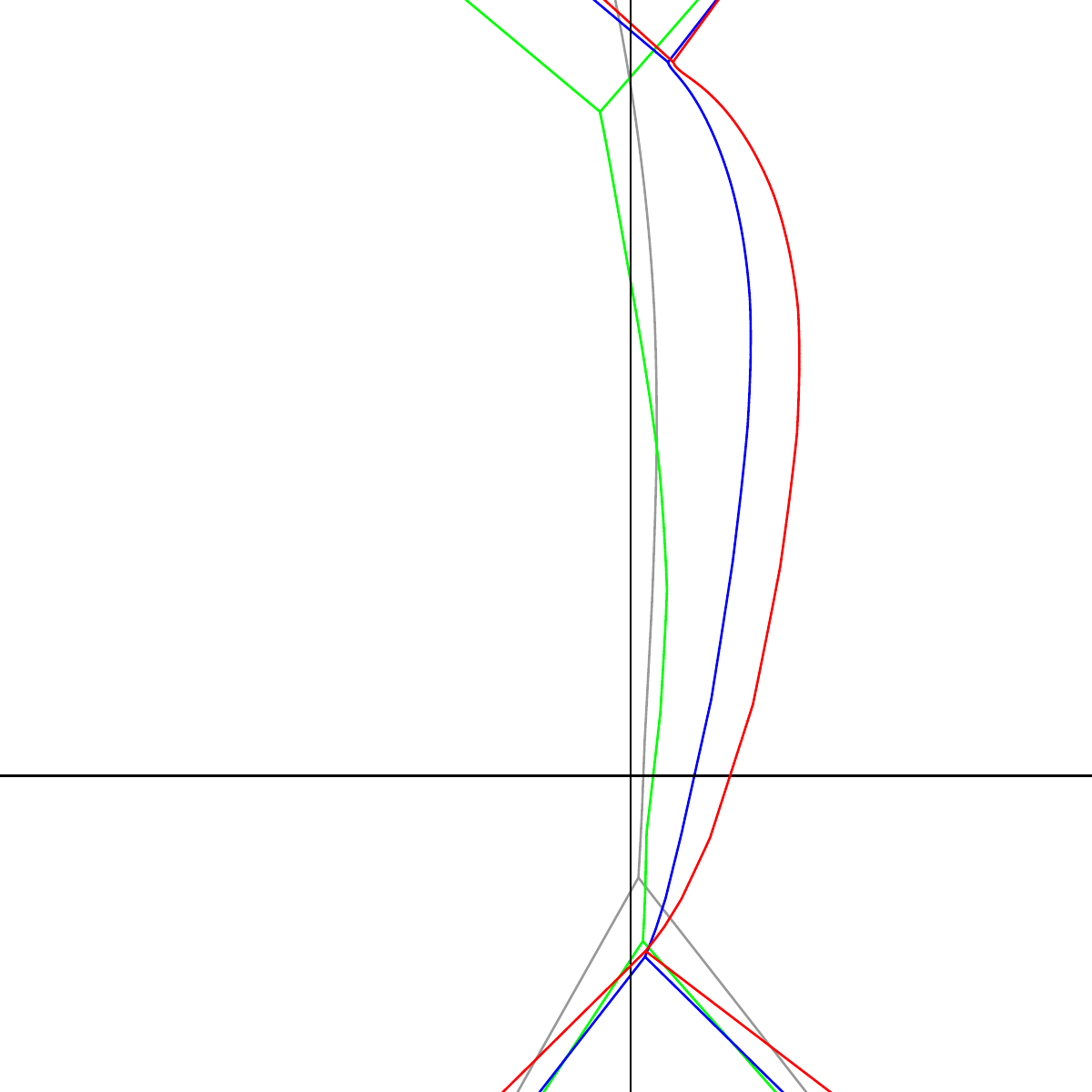}
        \caption{\cref{fig:triple-bubble-overlayed} zoomed in further reveals a small edge in all the bubbles.}
        \label{overlayed-edge}
    \end{subfigure}
    \hfill
    \begin{subfigure}[c]{0.4\textwidth}
        \centering
        \includegraphics[width=\linewidth]{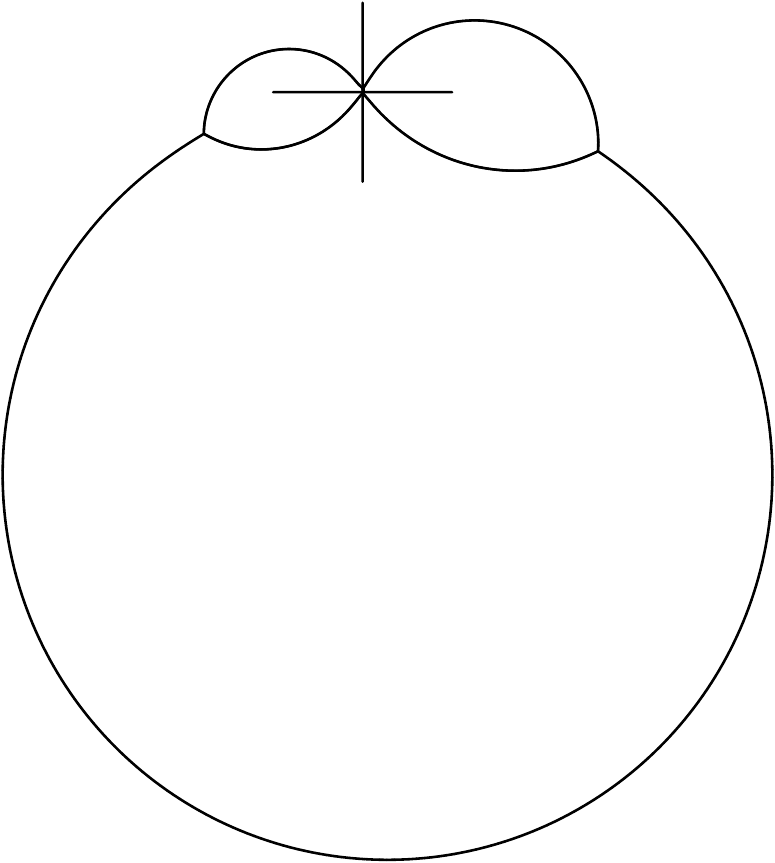}
        \caption{Volumes $0.1$, $0.5$, and $100$}
        \label{fig:triple-bubble-6}
    \end{subfigure}
    \hfill
    \begin{subfigure}[c]{0.2\textwidth}
        \centering
        \includegraphics[width=\linewidth]{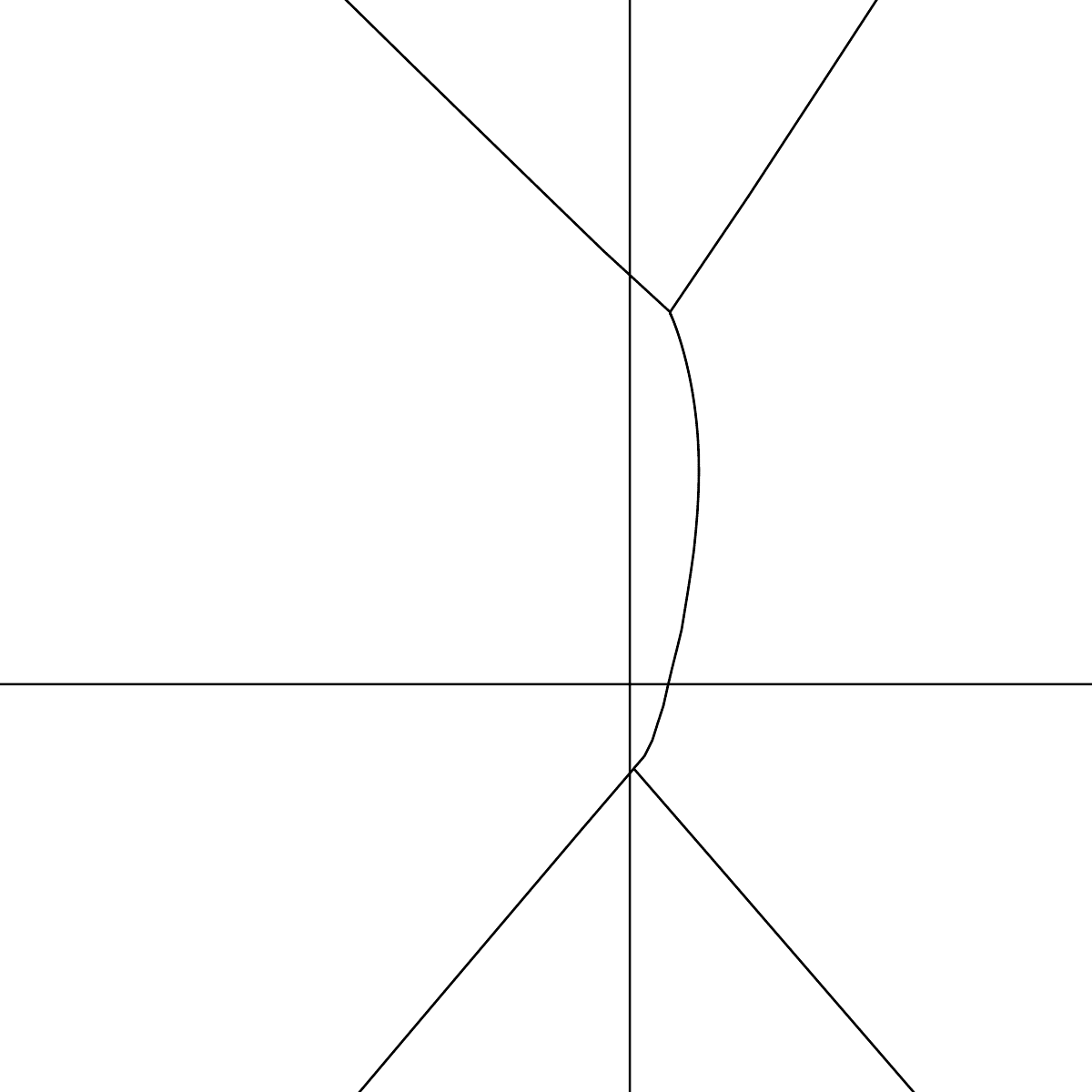}
        \caption{\cref{fig:triple-bubble-6} zoomed in on the origin}
        \label{fig:triple-bubble-8}
        \begin{minipage}{0.1cm}
        \vfill
        \end{minipage} 
    \end{subfigure}
    \caption{Computations in Brakke's Evolver \cite{Brakke2013} in $\R^2$ with density $r^2$ suggest that the optimal planar triple bubble consists of three circular arcs meeting at the origin and three nearly circular arcs meeting near the origin. Areas are labelled clockwise starting at the upper left bubble. Densities $r^3$, $r^4$, $r^5$, $r^6$, and $r^7$ are apparently similar.}
    \label{fig:triple-bubbles}
\end{figure}

\begin{figure}[h!]
    \centering
    \begin{subfigure}[c]{0.49\textwidth}
        \centering
        \includegraphics[width=.7\linewidth]{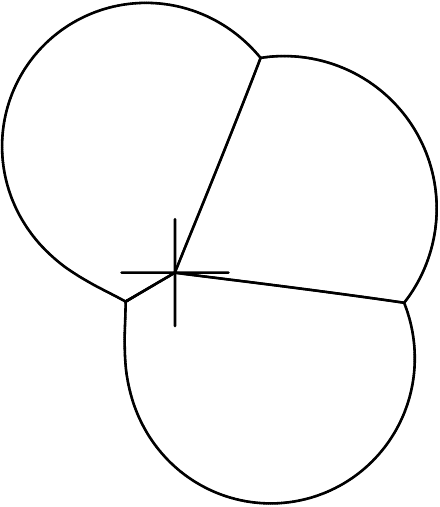}
        \caption{A triple bubble with equal areas for $p = 1.7$}
        \label{fig:1.7-triple}
    \end{subfigure}
    \hfill 
    \begin{subfigure}[c]{0.49\textwidth}
        \centering
        \includegraphics[width=.7\linewidth]{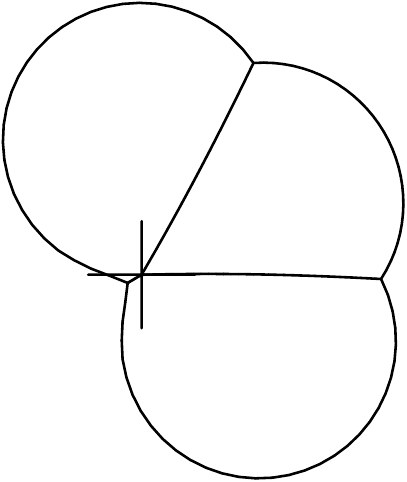}
        \caption{A triple bubble with equal areas for $p = 1.9$}
        \label{fig:1.9-triple}
        \begin{minipage}{0.1cm}
        \vfill
        \end{minipage} 
    \end{subfigure}
    \caption{As $p$ increases from $0$, one edge gets shorter, moving a second vertex near the origin. }
    \label{fig:triple-bubble-p<2}
\end{figure}

\begin{proposition}[Triple Bubble]
\label{tripleprop}
Computations with Brakke's Evolver \cite{Brakke2013} indicate that the optimal triple bubble in the plane with density $r^p$ consists of three circular arcs meeting at the origin, one shrinking as $p$ increases,  separating the bubbles from each other, and three constant-generalized-curvature curves (see \cref{genconstantcurvatureremark}), separating the bubbles from the exterior, as in \cref{fig:triple-bubbles} and \cref{fig:triple-bubble-p<2}.
\end{proposition}

\begin{figure}[h!]
    \centering
    \begin{subfigure}[c]{0.49\textwidth}
        \centering
        \includegraphics[width=0.45\linewidth]{Diagrams/10-10-10.pdf}
        \caption{Our conjectured triple bubble with equal volumes of $10$ has perimeter just over $63$.}
        \label{fig:triple-bubble-perimeter}
    \end{subfigure}
    \hfill 
    \begin{subfigure}[c]{0.49\textwidth}
        \centering
        \includegraphics[width=0.65\linewidth]{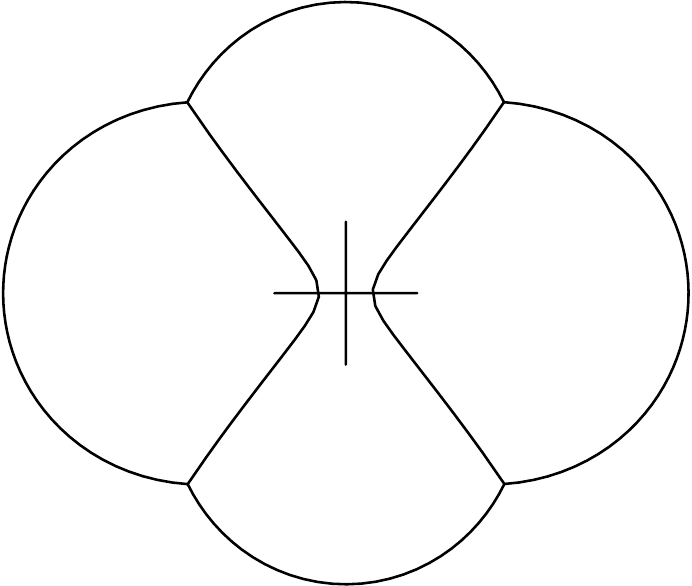}
        \caption{A linear chain with equal volumes of $10$ has perimeter just over $66$.}
        \label{fig:triple-chain-perimeter}
        \begin{minipage}{0.1cm}
        \vfill
        \end{minipage}
    \end{subfigure}
    \caption{Our conjectured triple bubble has less perimeter than a linear chain in the plane with density $r^2$. Densities $r^{0.5}$ and $r^3$ are apparently similar.}
    \label{fig:triple-comparison}
\end{figure}

\begin{figure}[h!]
    \centering
    \begin{subfigure}[c]{0.49\textwidth}
        \centering
        \includegraphics[width=.9\linewidth]{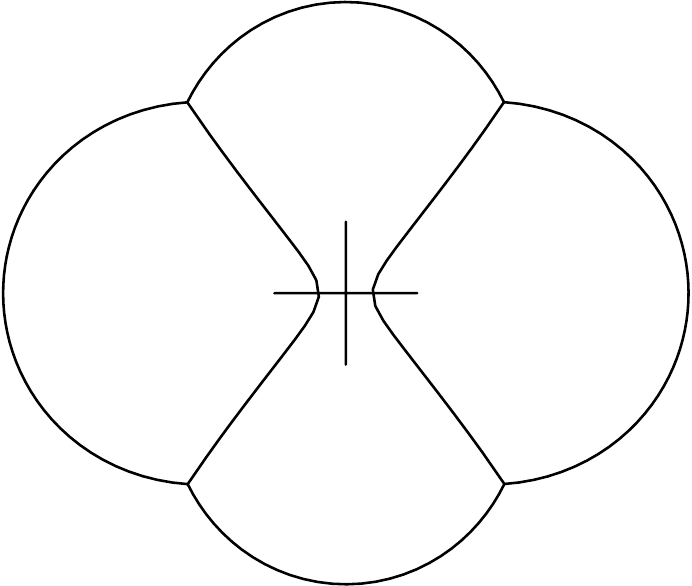}
        \label{starting-chain}
    \end{subfigure}
    \hfill
    \begin{subfigure}[c]{0.49\textwidth}
        \centering
        \includegraphics[width=.9\linewidth]{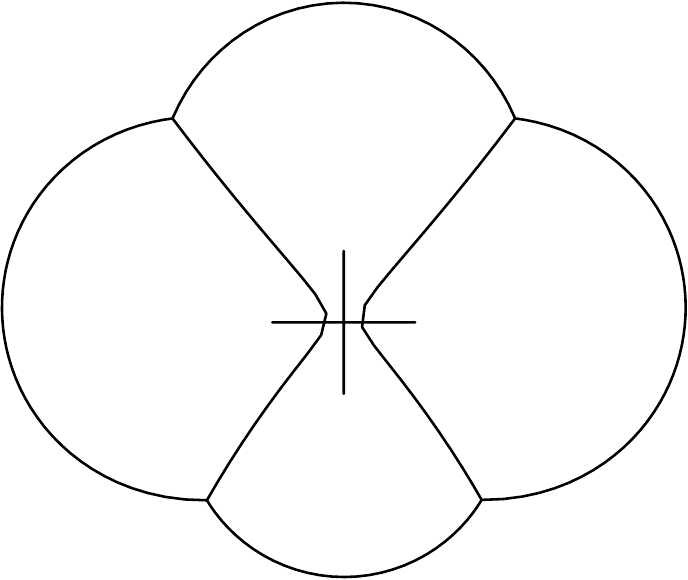}
        \label{fig:frame-1}
        \begin{minipage}{0.1cm}
        \vfill
        \end{minipage}
    \end{subfigure}
    \newline
    \begin{subfigure}[c]{0.49\textwidth}
        \centering
        \includegraphics[width=.9\linewidth]{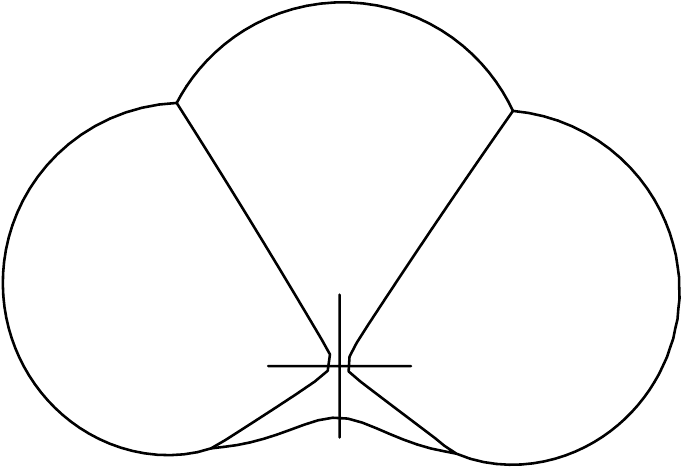}
        \label{fig:frame-2}
    \end{subfigure}
    \hfill
    \begin{subfigure}[c]{0.49\textwidth}
        \centering
        \includegraphics[width=.9\linewidth]{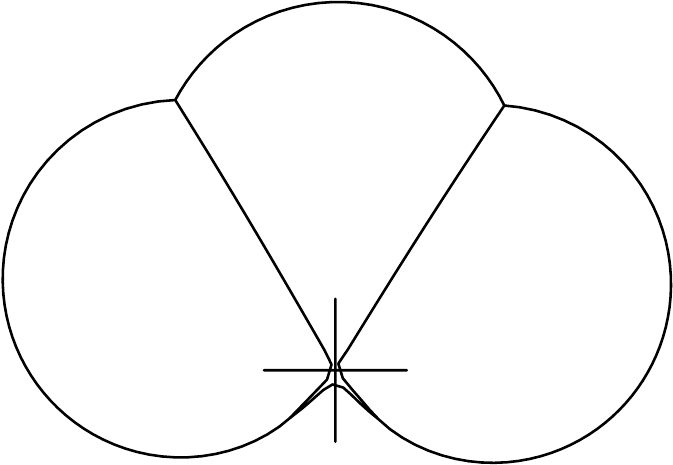}
        \label{fig:final-frame}
        \begin{minipage}{0.1cm}
        \vfill
        \end{minipage}
    \end{subfigure}
    \caption{A linear chain, after slight displacement upwards, evolves as far as possible (without changing its topological type) toward our conjectured triple bubble, here in the plane with density $r^2$. Densities $r^3$ and $r^4$ are apparently identical.}
    \label{fig:chain-evolution}
\end{figure}

\begin{proposition}
\label{triplechainprop}
Our conjectured triple bubble of \cref{tripleprop} has less perimeter than three bubbles in a linear chain, as in \cref{fig:triple-comparison}. The linear chain evolves toward our conjectured triple bubble as in \cref{fig:chain-evolution}.
\end{proposition}

\newpage

\begin{figure}[h]
    \centering
    \begin{subfigure}[c]{\textwidth}
        \centering         
        \includegraphics[width=0.49\linewidth]{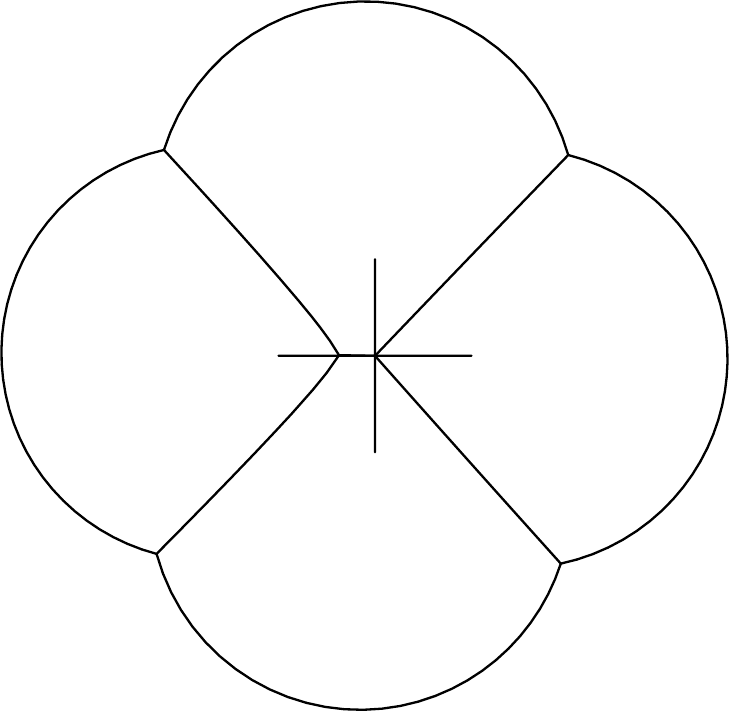}
        \caption{A quadruple bubble with equal areas for $p = 0.3$}
        \label{fig:0.5-quad}
    \end{subfigure}
    \newline
    \begin{subfigure}[c]{0.49\textwidth}
        \centering
        \includegraphics[width=\linewidth]{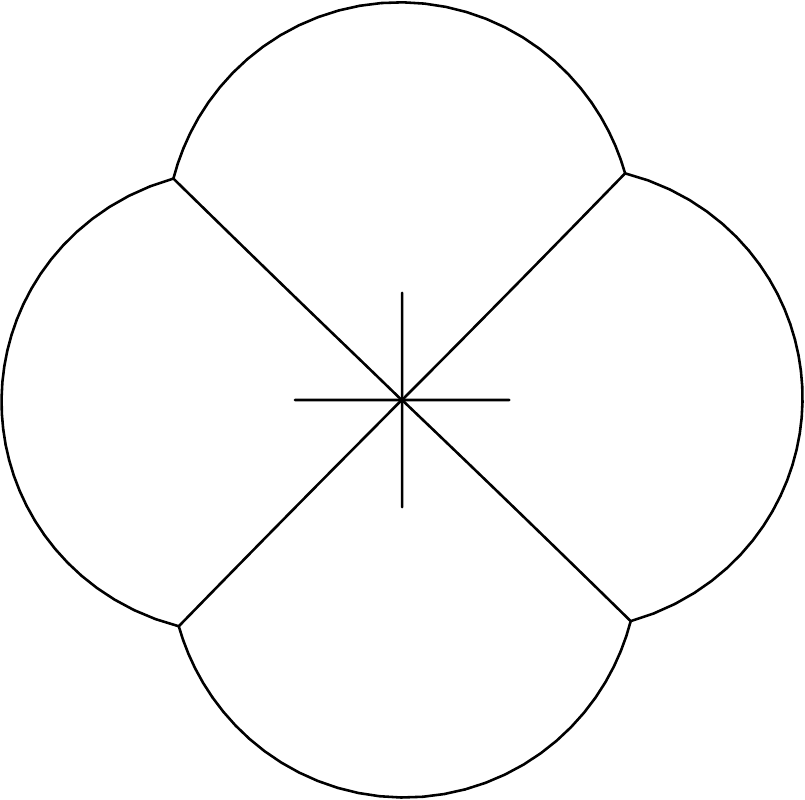}
        \caption{A quadruple bubble with equal areas for $p = 0.99$}
        \label{fig:0.99-quad}
    \end{subfigure}
    \hfill  
    \begin{subfigure}[c]{0.49\textwidth}
        \centering
        \hfill 
        \includegraphics[width=\linewidth]{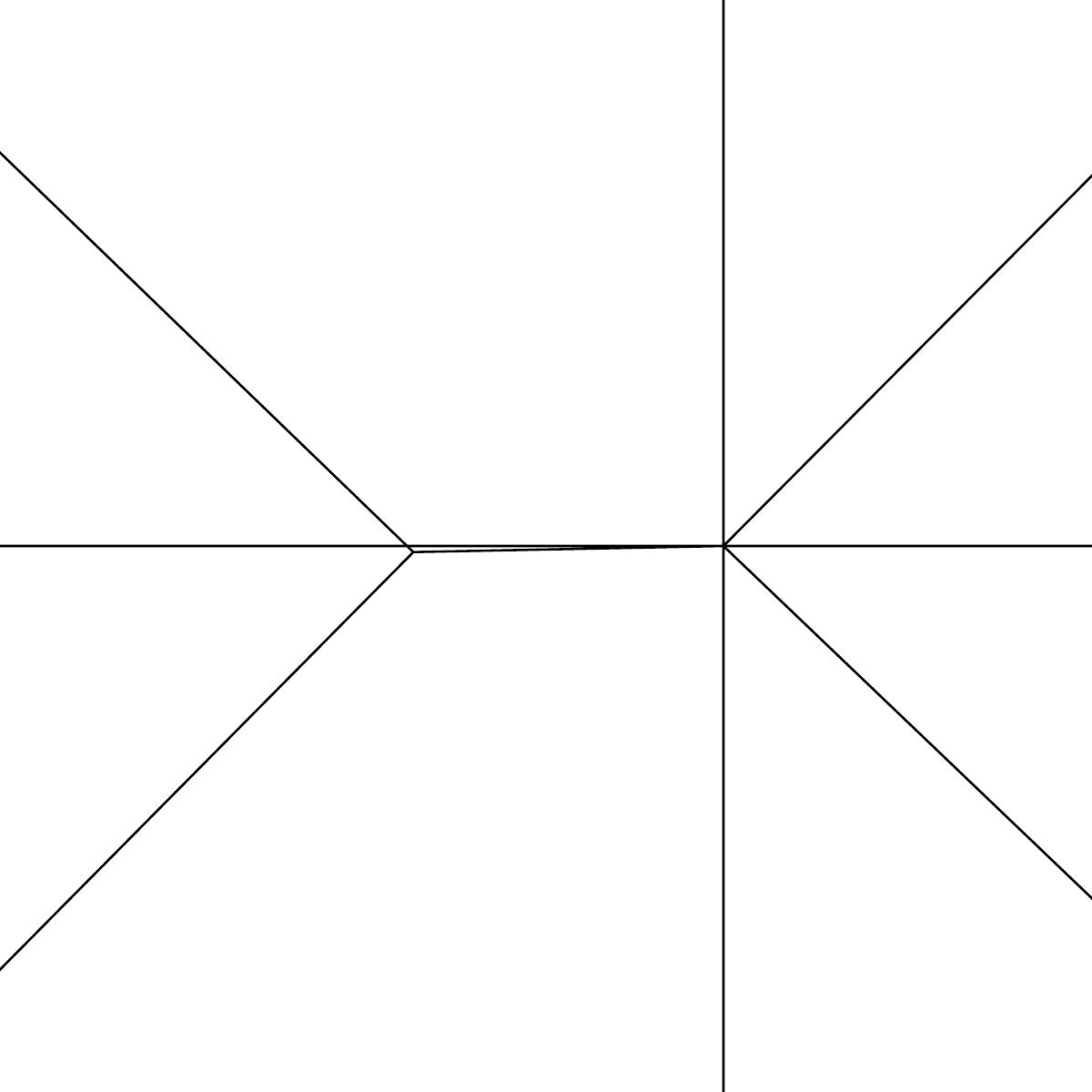}
        \caption{The same quadruple bubble as Figure \ref{fig:0.99-quad} zoomed in on the origin reveals a small edge}
        \label{fig:0.99-zoomed}
        \begin{minipage}{0.1cm}
        \vfill
        \end{minipage}
    \end{subfigure}
    \caption{Computations in Brakke's Evolver \cite{Brakke2013} in $\R^2$ with density $r^p$ for $p<1$ suggest that the optimal planar quadruple bubble has a short central edge with one endpoint at the origin, shrinking as $p$ increases toward $1$.}
    \label{fig:quadruple-bubbles-p<1}
\end{figure}

\begin{figure}[H]
    \centering
    \begin{subfigure}[c]{0.3\textwidth}
        \centering
        \includegraphics[width=4cm]{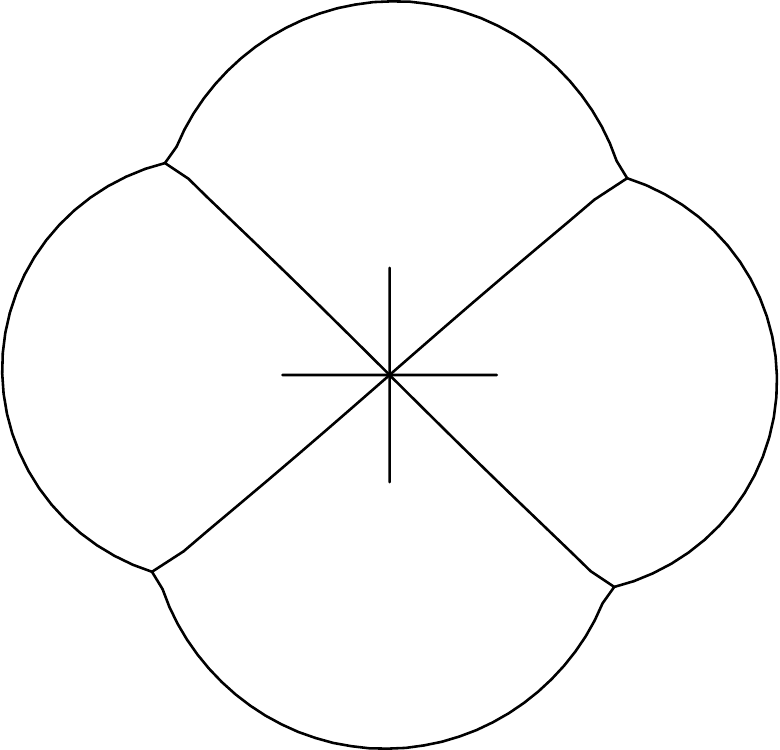}
        \label{equal_areas}
        \caption{Equal areas of $3$ show a symmetry absent in the Euclidean case (\ref{fig:euclidean-bubbles})}
    \end{subfigure}
    \hfill
    \begin{subfigure}[c]{0.3\textwidth}
        \centering
        \includegraphics[width=4cm]{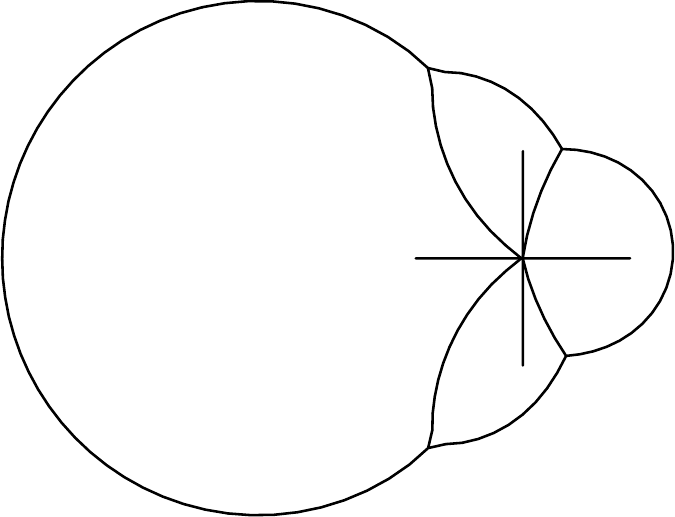}
        \label{fig:0.1-0.1-0.1-10}
        \caption{Areas of $0.1$, $0.1$, $0.1$, and $10$}
    \end{subfigure}
    \hfill 
    \begin{subfigure}[c]{0.3\textwidth}
        \centering
        \includegraphics[width=4cm]{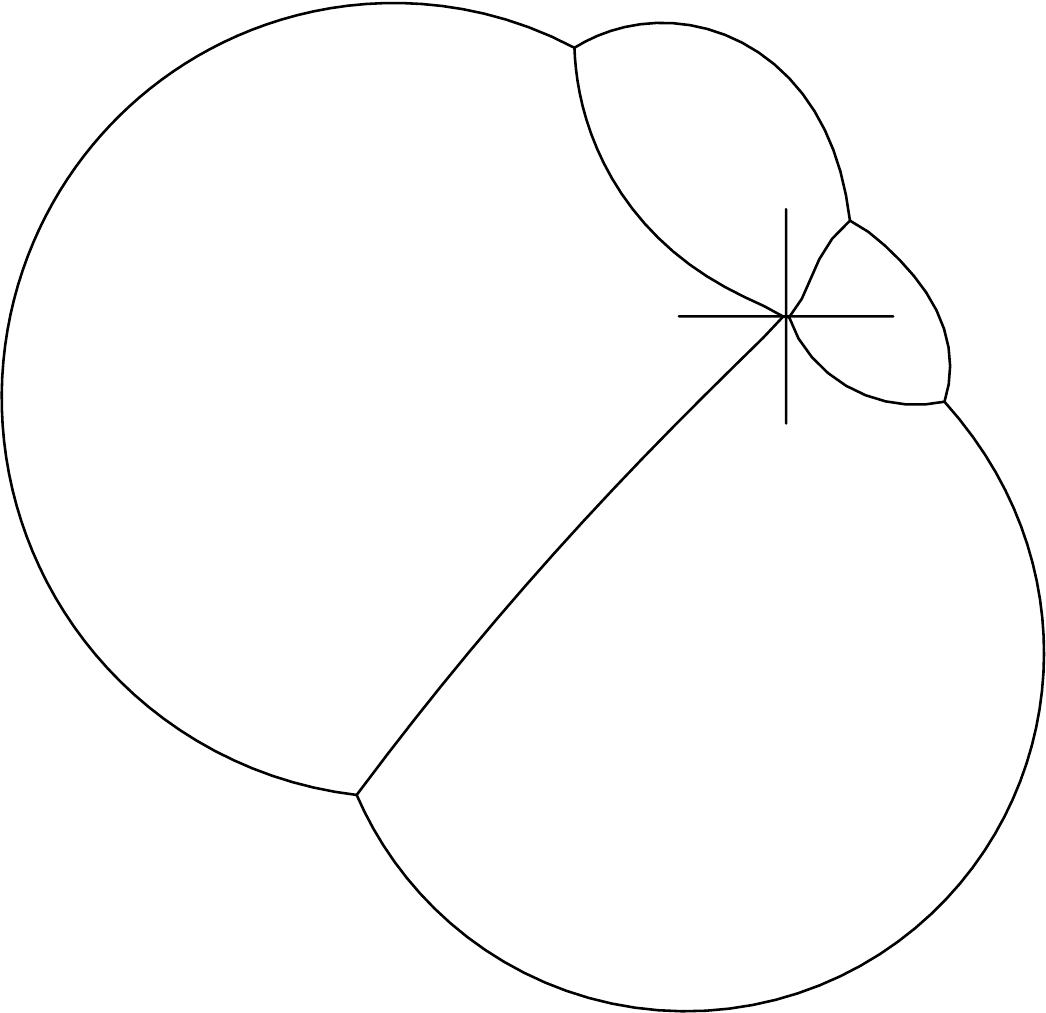}
        \label{fig:1-0.1-30-50}
        \caption{Areas of $1$, $0.1$, $30$, and $50$}
    \end{subfigure}
    \hfill
    \begin{subfigure}[c]{0.3\textwidth}
        \centering
        \includegraphics[width=4cm]{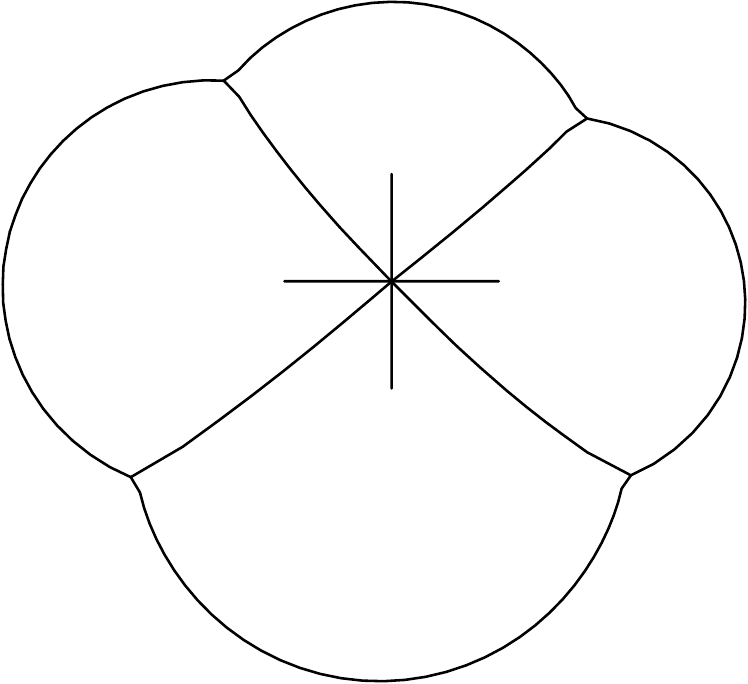}
        \label{fig:1-2-4-3}
        \caption{Areas of $1$, $2$, $4$, and $3$}
    \end{subfigure}
    \hfill 
    \begin{subfigure}[c]{0.3\textwidth}
        \centering
        \includegraphics[width=4cm]{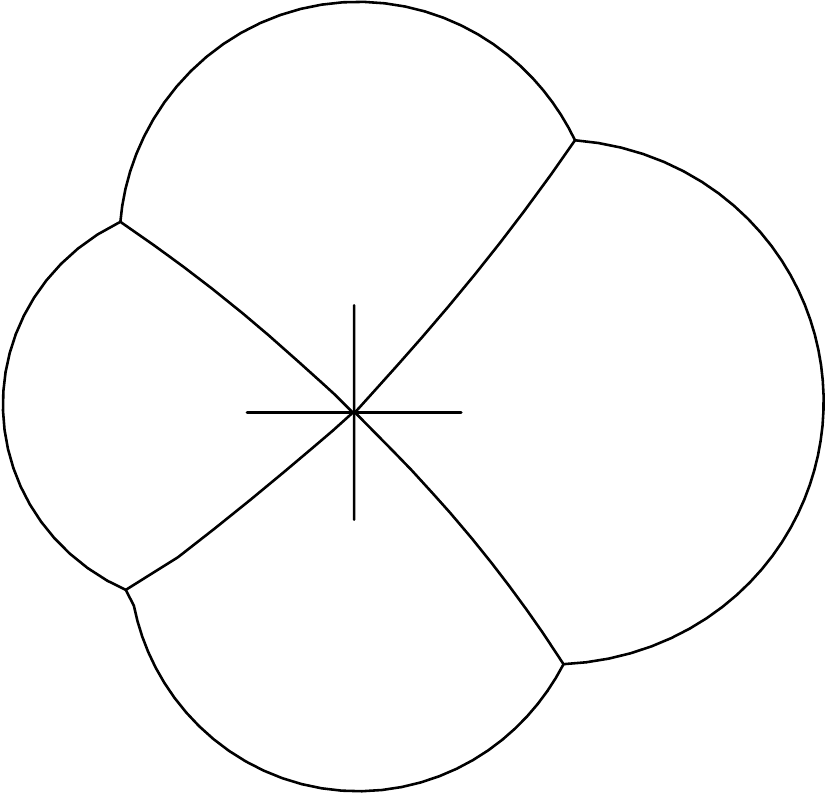}
        \label{fig:4-7-3-2}
        \caption{Areas of $4$, $7$, $3$, and $2$}
    \end{subfigure}
    \hfill
    \begin{subfigure}[c]{0.3\textwidth}
        \centering
        \includegraphics[width=4cm]{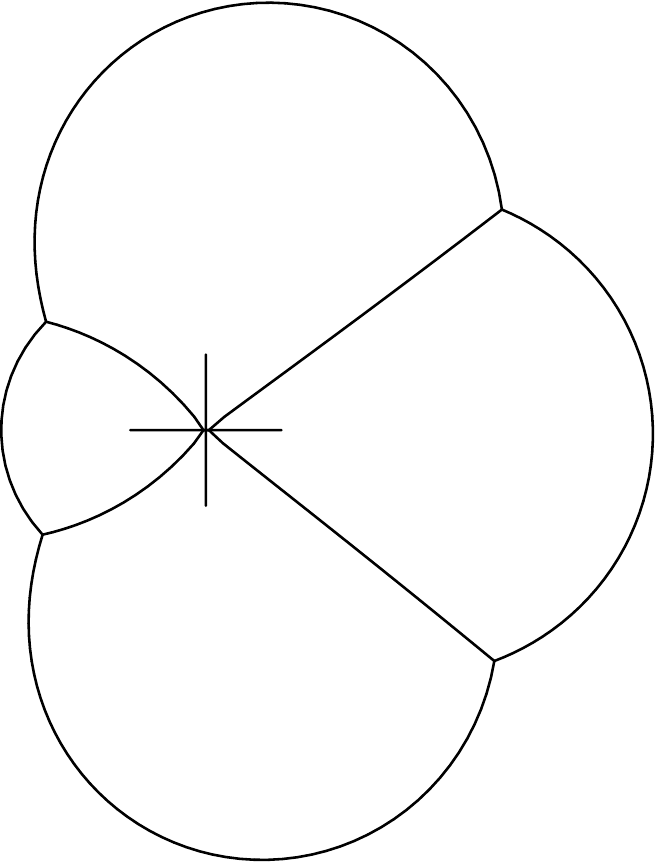}
        \caption{Areas of $20$, $20$, $20$ and $1$}
            \label{fig:20-20-20-1}

    \end{subfigure}
    \newline
    \begin{subfigure}[c]{0.3\textwidth}
        \centering
        \includegraphics[width=4cm]{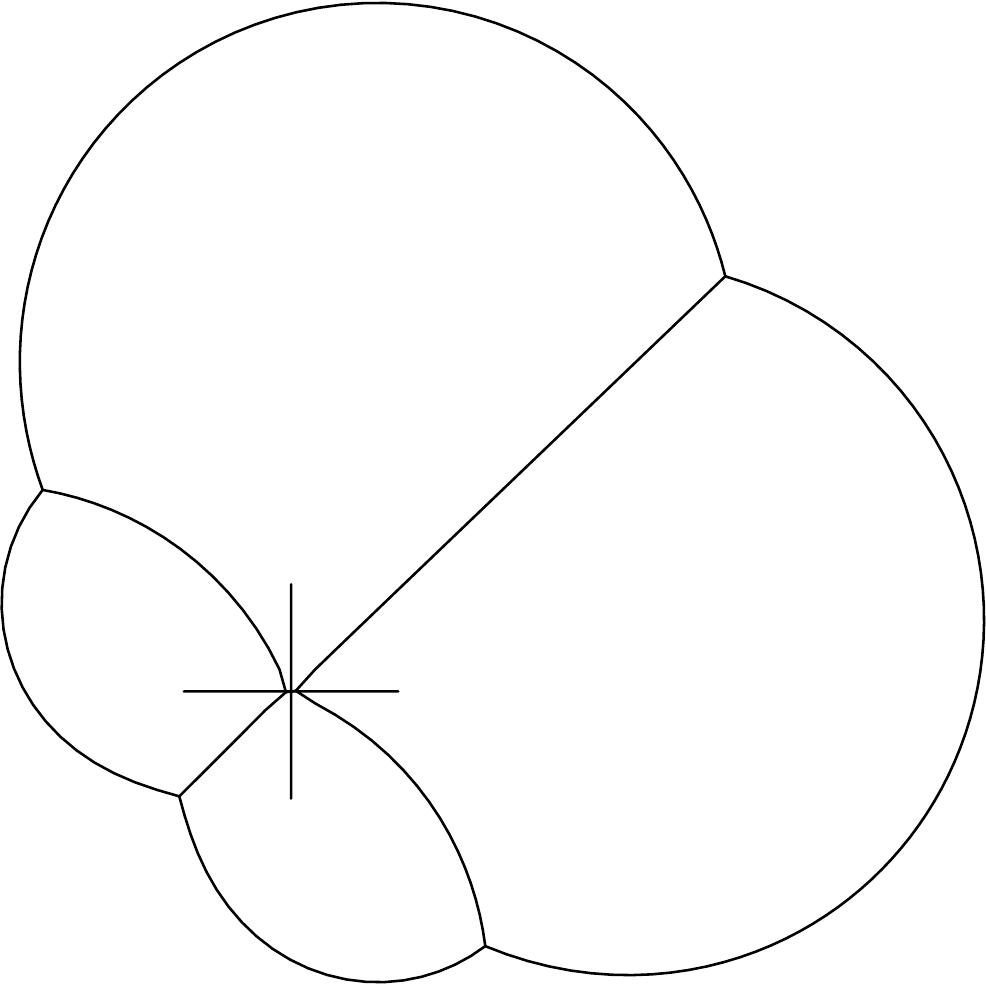}
        \caption{Areas of $30$, $30$, $1$, and $1$, with total perimeter just under $104$}
        \label{fig:30-30-1-1}
    \end{subfigure}
    \hfill
    \begin{subfigure}[c]{0.3\textwidth}
        \centering
        \includegraphics[width=3cm]{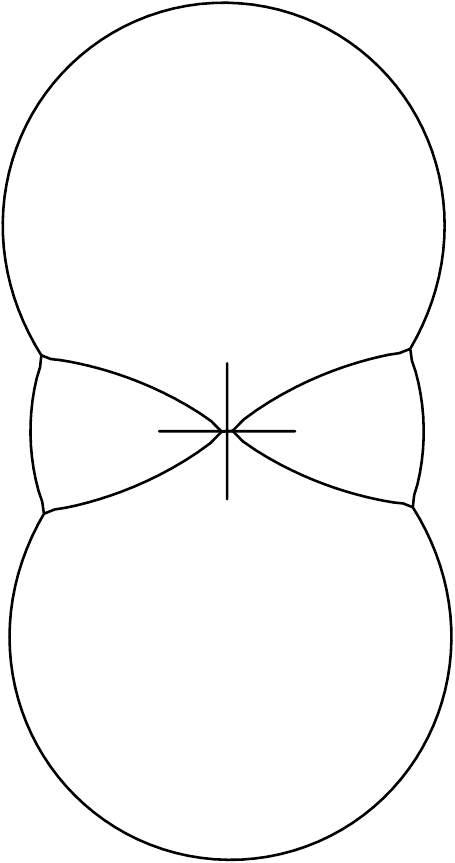}
        \label{fig:30-1-30-1}
        \caption{An inferior version of \cref{fig:30-30-1-1} with areas of $30$, $1$, $30$, and $1$, with total perimeter just over $106$}
    \end{subfigure}
    \hfill
    \begin{subfigure}[c]{0.3\textwidth}
        \centering
        \includegraphics[width=4cm]{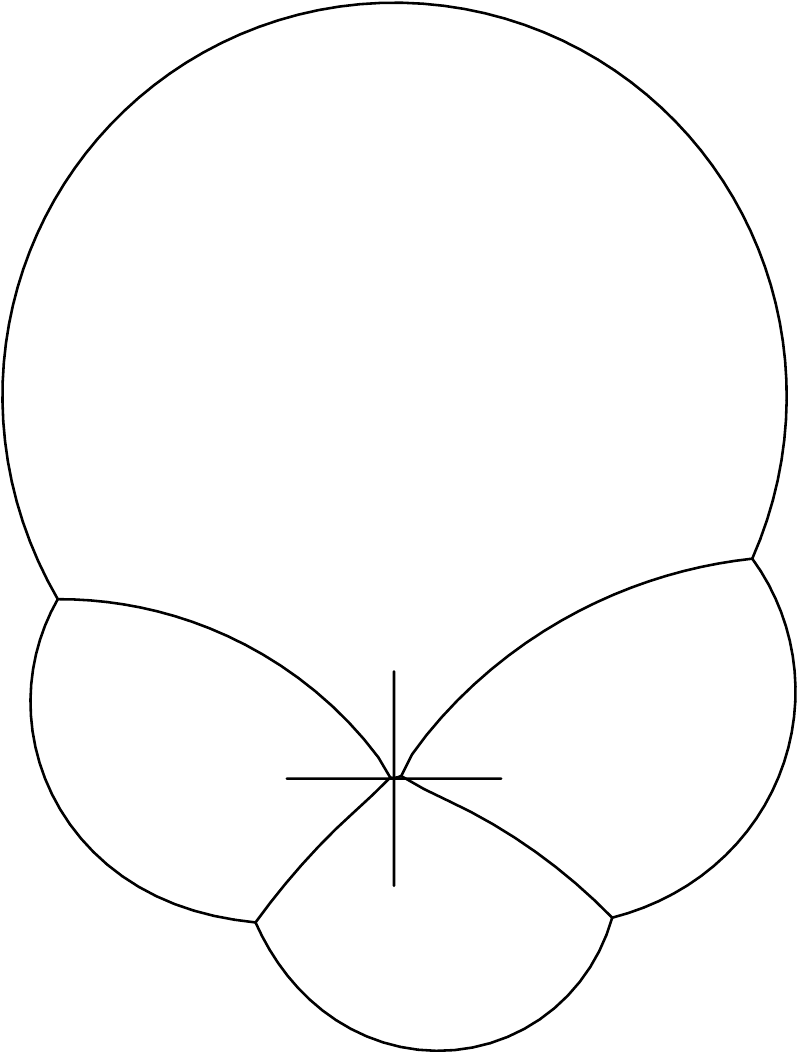}
        \label{fig:50-3-1-2}
        \caption{Areas of $50$, $3$, $1$, and $2$}
    \end{subfigure}
    \caption{For $p \ge 1$, the central edge has collapsed and four circular arcs meet at the origin. The areas are labelled clockwise, starting from the top bubble. Densities $r^3$ and $r^4$ are apparently identical.}
    \label{fig:quad-bubbles}
\end{figure}

\begin{proposition}[Quadruple Bubble]
\label{quadprop}
Computations with Brakke's Evolver \cite{Brakke2013} indicate that for the optimal quadruple bubble in the plane with density $r^p$ as $p$ increases from $0$ (the standard Euclidean quadruple bubble as in \cref{fig:euclidean-bubbles}) the central edge with one endpoint at the origin shrinks as in \cref{fig:quadruple-bubbles-p<1} until it disappears when p reaches 1, after which four circular arcs meet at the origin (where the density vanishes) as in \cref{fig:quad-bubbles}. 
\end{proposition}

\begin{figure}[h!]
    \centering
    \begin{subfigure}[c]{0.49\textwidth}
        \centering
        \includegraphics[width=\linewidth]{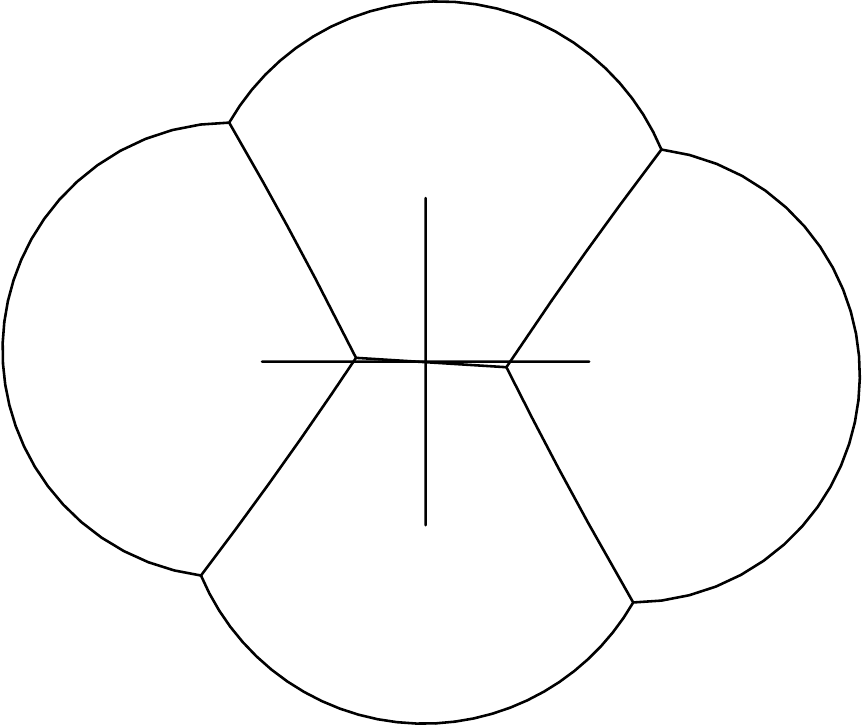}
        \caption{The standard Euclidean quadruple bubble with unit Euclidean areas,  centered at the origin in the plane with density $r^2$, has areas around $0.5, 0.8, 0.5$ and $0.8$ and perimeter around $10.86$}
        \label{fig:euclidean-quadruple-perimeter}
    \end{subfigure}
    \hfill 
    \begin{subfigure}[c]{0.49\textwidth}
        \centering
        \includegraphics[width=\linewidth]{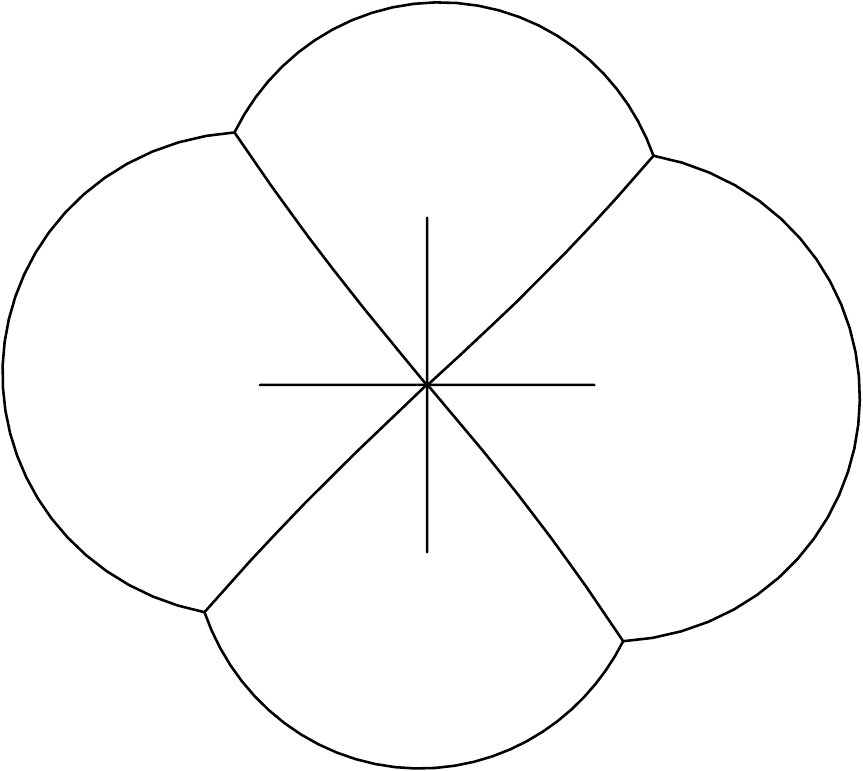}
        \caption{Our conjectured quadruple bubble with the same areas in the plane with density $r^2$ has perimeter around $10.81$}
        \label{fig:conjectured-quadruple-perimeter}
        \begin{minipage}{0.1cm}
        \vfill
        \end{minipage}
    \end{subfigure}
    \caption{With the same areas in the plane with density $r^2$, our conjectured quadruple bubble has less perimeter than the standard Euclidean quadruple bubble.}
    \label{fig:quadruple-comparison}
\end{figure}

For example, in the plane with density $r^p$ for $p \ge 1$ our conjectured quadruple bubble of \cref{quadprop} has less perimeter than the Euclidean quadruple bubble, as in \cref{fig:quadruple-comparison}.

\newpage 

\begin{figure}[H]
    \centering
    \begin{subfigure}[c]{\textwidth}
        \centering
        \includegraphics[width=0.3\linewidth]{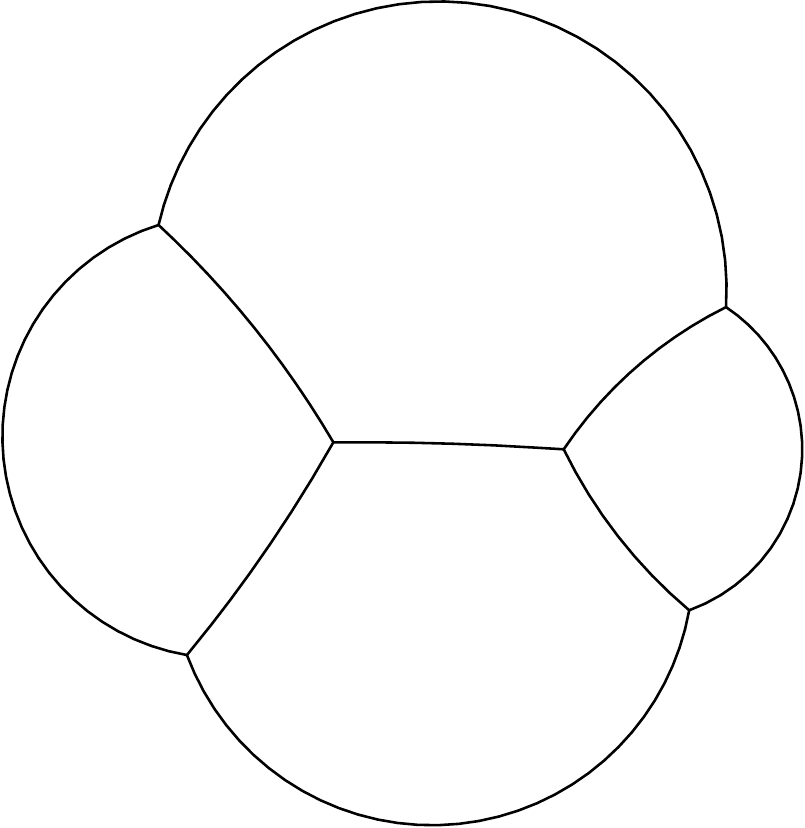}
        \caption{The optimal Euclidean quadruple bubble apparently has the two largest regions on opposite sides of the central edge.}
        \label{fig:euclidean-ordering}
    \end{subfigure}
    \newline    
    \begin{subfigure}[c]{0.49\textwidth}
        \centering
        \includegraphics[width=0.75\linewidth]{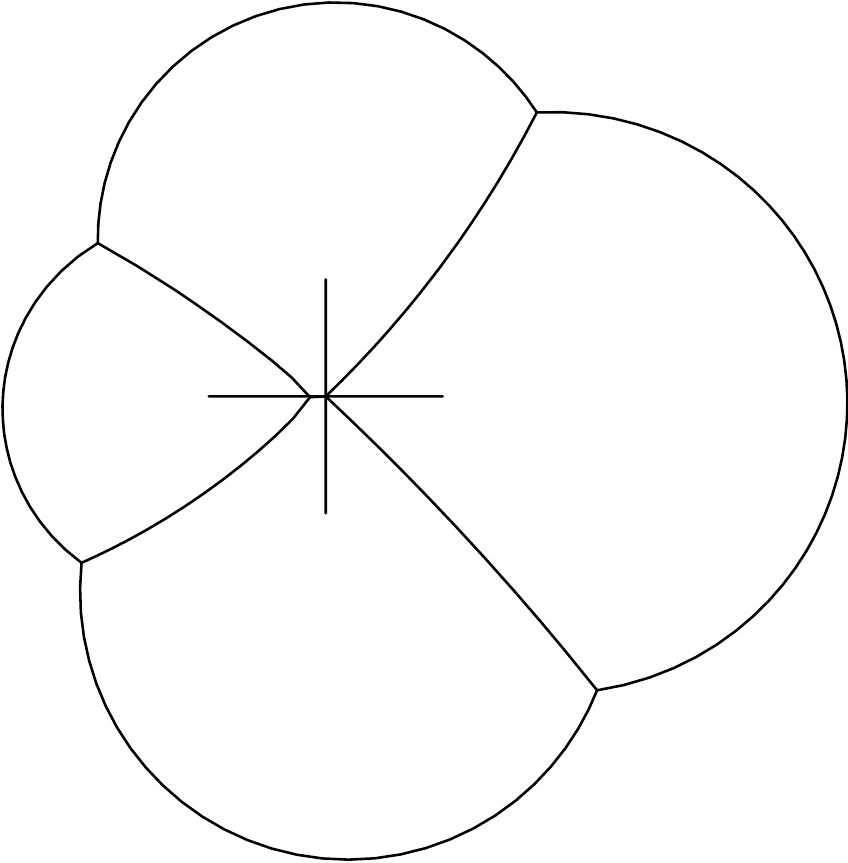} 
        \caption{As $p$ increases towards $1$, the optimal quadruple bubble has the largest and smallest regions on opposite ends of the central edge.}
        \label{fig:0.5-ordering}
    \end{subfigure} 
    \hfill
    \begin{subfigure}[c]{0.49\textwidth}
        \centering
        \includegraphics[width=0.75\linewidth]{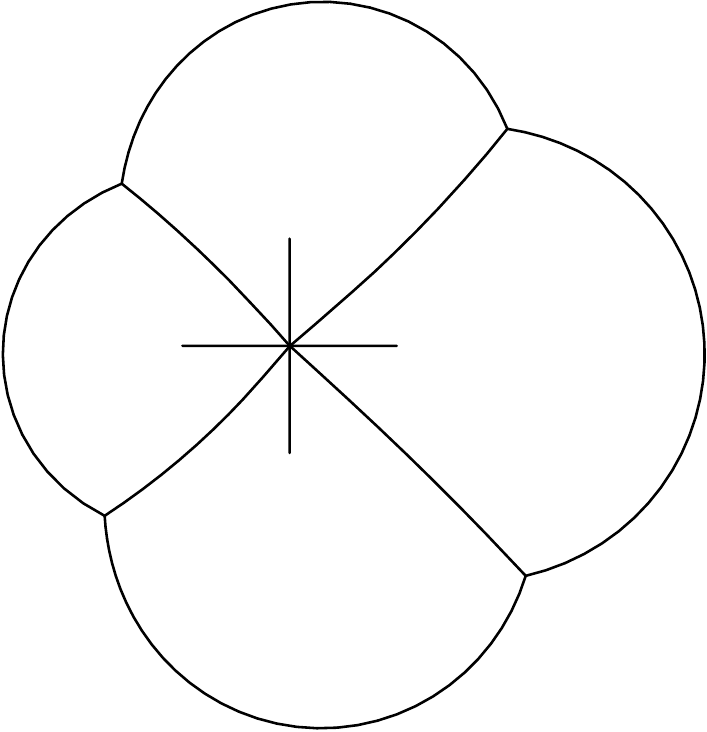}
        \caption{Once the central edge has disappeared, for $p \ge 1$, the largest and smallest bubbles remain opposite.}
        \label{fig:2-ordering}
        \begin{minipage}{0.1cm}
        \vfill
        \end{minipage}
    \end{subfigure}
    \caption{Optimal quadruple bubble orderings for different values of $p$.}
    \label{fig:quadruple-arrangements}
\end{figure}

\begin{proposition}
\label{quadorientationprop}
Computations indicate that the optimal Euclidean quadruple bubble has the two largest areas on opposite sides of the central edge.  As $p$ increases towards $1$ it prefers the largest and smallest bubbles on opposite ends of the central edge. Once the central edge has disappeared, for $p \ge 1$, the largest and smallest bubbles remain opposite. See \cref{fig:quadruple-arrangements}.
\end{proposition}

This arrangement of unequal areas is apparently a new conjecture even for the Euclidean case. It indicates that that the Euclidean quadruple bubble never has convex regions, so the results of \textcite{DeRosa2023} never apply. It is known to be true for small deformations of the equal-areas minimizer because it saves perimeter to shrink convex regions and expand nonconvex regions.

\begin{figure}[h!]
    \centering
    \begin{subfigure}[c]{0.48\textwidth}
        \centering
        \includegraphics[width=0.49\linewidth]{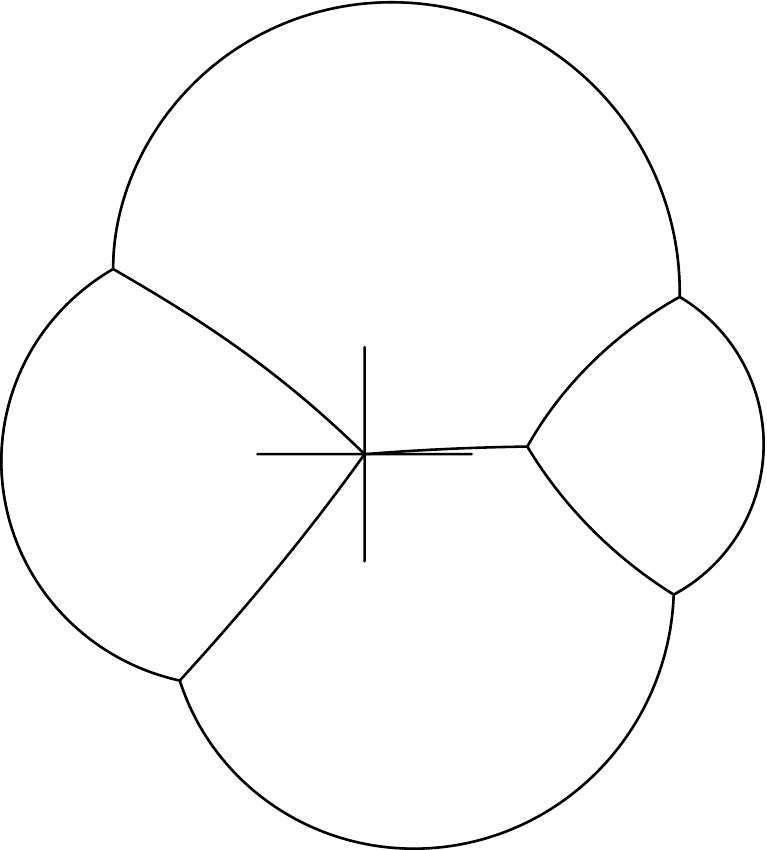}
        \caption{Our conjectured quadruple bubble with the vertex of the larger bubbles on the origin has perimeter around $17.67$.}
        \label{fig:origin-big-bubbles}
    \end{subfigure}
    \hfill
    \begin{subfigure}[c]{0.48\textwidth}
        \centering
        \includegraphics[width=0.49\linewidth]{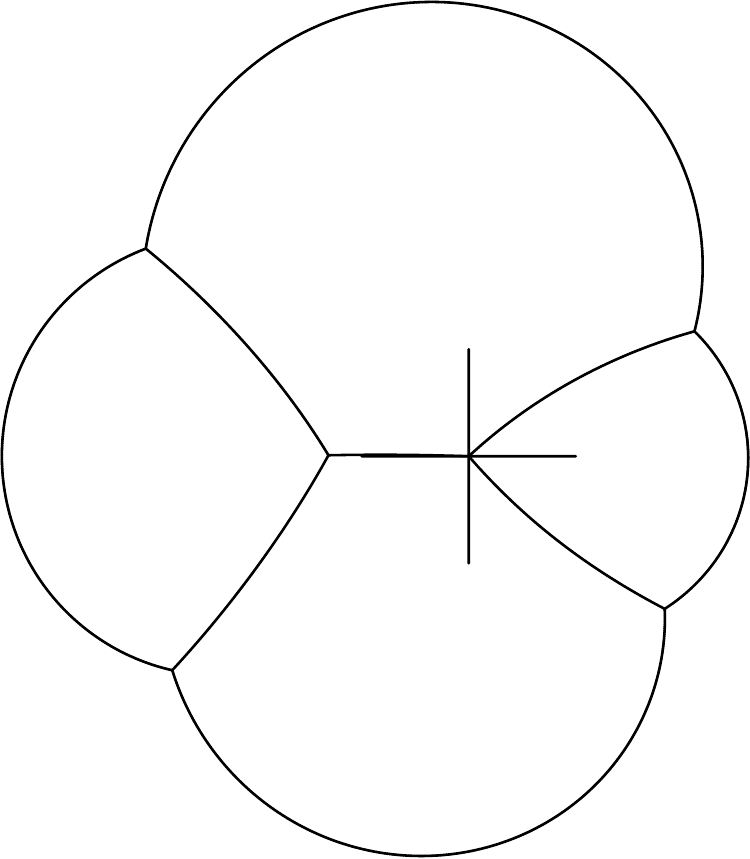}
        \caption{The quadruple bubble with the vertex of the smallest bubble on the origin has perimeter around $17.70$.}
        \label{fig:origin-small-bubbles}
        \begin{minipage}{0.2cm}
        \vfill
        \end{minipage}
    \end{subfigure}
    \caption{Our conjectured quadruple bubble has the vertex of the larger bubbles on the origin (here for density $r^{0.1}$).}
    \label{fig:quadruple-origin}
\end{figure}

\begin{proposition}
\label{quadoriginprop}
Computations indicate that when a central edge is present, our conjectured quadruple bubble of \cref{quadprop} is most effective with the vertex of the larger bubbles on the origin, as in \cref{fig:quadruple-origin}.
\end{proposition}

This is the opposite of the $\R^1$ case where the vertex of the smallest bubbles is on the origin (Ross \cite[Thm.~1]{Ross2022}). In $\R^2$ larger bubbles require more perimeter, so it makes sense to have them closer to the origin. 

\clearpage 
\newpage
\printbibliography
\end{document}